\documentclass[11pt,oneside,a4paper]{article}

\usepackage[french]{babel}
\usepackage[latin1]{inputenc}
\usepackage{amsmath}
\usepackage{amsthm}
\usepackage{amssymb}
\usepackage{array}
\usepackage[dvips]{graphicx}
\newtheorem{thm}{Th\'eor\`eme}[section]
\newtheorem{prop}[thm]{Proposition} 
\newtheorem{lem}[thm]{Lemme}
\newtheorem{cor}[thm]{Corollaire}
 
\theoremstyle{definition}
 
\newtheorem{defn}[thm]{D\'efinition}

\theoremstyle{remark}
 
\newtheorem{rem}[thm]{Remarque}

\def\z{{\widehat z }}
\def\tpsi{{\widetilde \Psi }}

\title{R\'esurgence des solutions BKW formelles d'une EDO
  singuli\`erement perturb\'ee}
\author{Jean-Marc Rasoamanana}
\date{D\'epartement de Math\'ematiques, UMR CNRS 6093,
Universit\'e d'Angers, 2 Boulevard Lavoisier, 49045 Angers Cedex 01,
France.}
 
\begin{document}
 \maketitle 
 \newpage
 \tableofcontents
 \newpage

\section{Introduction}

\subsection{Pr\'esentation}

Les EDO singuli\`erement perturb\'ees servent tr\`es souvent de
mod\`eles, notamment en physique quantique (le param\`etre de
perturbation $\varepsilon$ repr\'esentant alors la constante de Planck
$\hbar$ des physiciens).\\
Un exemple classique est l'\'equation de Schr\"odinger
unidimensionnelle stationnaire dans le champ complexe~:
\begin{equation}\label{eq1intro}
\varepsilon^2 \frac{d^2 Y}{d q^2} -V(q) Y =0,
\end{equation}
o\`u la fonction potentielle $V$ est analytique (par exemple
polynomiale).

L'\'etude de telles \'equations conduit de mani\`ere naturelle \`a
consid\'erer des solutions formelles (en $\varepsilon$) qu'on appelle
d\'eveloppements BKW (du nom des physiciens Brillouin, Krammers et
Wentzel) ou d\'eveloppements semi-classiques.

D'une mani\`ere g\'en\'erale, ces d\'eveloppements formels sont
divergents, ce qui conduit alors \`a \'etudier leur caract\`ere
r\'esurgent ou sommable (de Borel) par rapport au param\`etre de
perturbation $\varepsilon$ (ce qu'\'Ecalle appelle r\'esurgence
quantique ou co\'equationnelle dans \cite{Ec85}).

Les techniques de sommation ont \'et\'e largement d\'evelopp\'ees,
notamment gr\^ace aux travaux de J.P. Ramis (\cite{R78}, \cite{R85} et
\cite{R84} notamment) et de J. Ecalle
(\cite{Ec81-1}, \cite{Ec81-2}, \cite{Ec84} et \cite{Ec85} par exemple), et utilis\'ees avec succ\`es pour retrouver, \`a partir de
certains d\'eveloppements formels, des ``vraies'' solutions exactes de
l'\'equation consid\'er\'ee. De fait, l'int\'er\^et de la sommation de
Borel, notamment dans la m\'ethode BKW, est immense, tant au niveau
math\'ematique proprement dit (voir \cite{DP99}, \cite{Din73}, \cite{Ec84} ou \cite{Vo83}) qu'au niveau des
applications en physique (voir \cite{BB74}, \cite{BebWu}, \cite{Sim82} et \cite{ZJ89}  par exemple).

Cette m\'ethode de sommation dans le cadre BKW est souvent qualifi\'ee
d'analyse BKW exacte (ou d'analyse semi-classique exacte) et cette
exactitude permet notamment l'obtention, dans certains cas, de
formules de connexion entre les diff\'erentes solutions BKW (voir
\cite{Vo83} par exemple).

En ce qui concerne l'aspect r\'esurgent de tels d\'eveloppements, il
appara\^it que les solutions BKW peuvent \^etre per\c cues comme un
v\'eritable codage exact de vraies solutions (voir \cite{DDP97}) et non pas
seulement comme de simples approximations. Le ph\'enom\`ene de
Stokes s'interpr\`ete alors naturellement comme discontinuit\'e dans
de tels codages.

Un th\'eor\`eme d'Ecalle affirme que dans le cas de l'\'equation
(\ref{eq1intro}), il existe toujours une base de solutions BKW formelles
r\'esurgentes, pourvu que le potentiel $V$ se comporte ``suffisamment
bien \`a l'infini''. Toutefois, de l'avis des sp\'ecialistes, ce
th\'eor\`eme n'est pas encore compl\`etement d\'emontr\'e dans sa
g\'en\'eralit\'e.

Notre point de vue s'inscrit dans ce ``courant de pens\'ee''.\\
Le sujet principal de cet article est l'\'etude de l'\'equation
diff\'erentielle ordinaire singuli\`erement perturb\'ee~:
\begin{equation}\label{eq4intro}
\displaystyle  \frac{d^2 \Phi}{d z^2} - \frac{z}{\varepsilon^2} \Phi = F(z) \Phi,
\end{equation} 
o\`u $F$ d\'esigne une fonction holomorphe, au moins au voisinage de
l'origine, et $\varepsilon$ est un petit param\`etre complexe.\\
Remarquons que cette \'equation ne rentre pas dans le champ
d'applications du th\'eor\`eme d'Ecalle
pr\'ec\'edemment cit\'e.

En utilisant les outils de la th\'eorie BKW
exacte, nous allons analyser les propri\'et\'es de r\'esurgence
(param\'etrique) d'une classe de solutions BKW formelles ``bien
normalis\'ees''.\\
Nous discuterons \'egalement de leur \'eventuel
caract\`ere sommable.\\
En outre, le r\'esultat principal de cet article est le th\'eor\`eme
suivant~:
\begin{thm}
Lorsque $F$ est holomorphe au voisinage de l'origine (respectivement enti\`ere), il existe une
famille de solutions BKW \'el\'ementaires r\'esurgentes de type Airy
local (respectivement de type Airy) $\Phi_{bkw}(z,\varepsilon)$ de l'\'equation (\ref{eq4intro}),
\end{thm}
au sens de la d\'efinition suivante~:
\begin{defn}\label{deftypeAiry}
Un symbole r\'esurgent \'el\'ementaire $\Phi(z,\varepsilon)$ est dit
de type Airy local en $(z=0,\varepsilon=0)$ (respectivement de type Airy)
s'il v\'erifie les conditions suivantes~:
\begin{enumerate}
\item son support singulier est inclus dans la courbe alg\'ebrique
  \linebreak $\mathcal{C} = \{(z, \xi),\,\, 9\xi^2 = 4z^3 \}$ 
au voisinage de $(z=0,\varepsilon=0)$,\\
\item pour toute direction $\alpha$, et tout germe de secteur de
  Stokes (respectivement secteur de Stokes) $S$ relatif \`a $\alpha$, toute d\'etermination du symbole
  $\Phi(z,\varepsilon)$ peut s'\'ecrire comme la d\'ecomposition
  locale (respectivement d\'ecomposition), pour $z \in S$, d'une fonction confluente (respectivement
  d'une fonction confluente r\'esurgente) \`a support
  singulier inclus dans $\mathcal{C}$.
\end{enumerate}
\end{defn}

Par ailleurs, l'une de nos motivations est d'appliquer nos r\'esultats
\`a l'\'equation de Schr\"odinger (\ref{eq1intro}) : en effet, cette
derni\`ere se ram\`ene \`a notre \'equation principale
(\ref{eq4intro}) via un changement de variable analytique.\\
En particulier, nous \'etablissons un th\'eor\`eme local r\'esurgent
de r\'eduction (au voisinage d'un point tournant simple) qui affirme
que l'\'equation (\ref{eq1intro}) peut se ramener \`a l'\'equation
d'Airy~:
\begin{equation}\label{eq3intro}
\frac{d^2 y}{d s^2} =\frac{s}{\varepsilon^2} y,
\end{equation}
(i.e l'\'equation d'Airy est le mod\`ele local universel pour un point
tournant simple).

\subsection{Contenu}

Le papier est organis\'e de la mani\`ere suivante.\\
Dans un premier temps, nous allons analyser en d\'etail dans la section \ref{sec2}  l'\'equation
(\ref{eq4intro}) dans le cas o\`u $F=0$. : l'\'equation (\ref{eq4intro})
n'est alors rien d'autre que l'\'equation d'Airy, qui va nous servir
de mod\`ele pour l'analyse BKW exacte de l'\'equation (\ref{eq4intro})
dans le cas g\'en\'eral. En particulier, nous y d\'efinissons le
symbole BKW d'Airy, y rappelons ses propri\'et\'es de r\'esurgence
et sommabilit\'e et analysons en d\'etail le ph\'enom\`ene de Stokes associ\'e.

Dans la section \ref{sec3}, nous commen\c cons par l'analyse BKW formelle de
l'\'equation (\ref{eq4intro}) dans le cas g\'en\'eral en montrant
l'existence d'une famille de solutions BKW formelles "bien
normalis\'ees" de (\ref{eq4intro}).

La section \ref{sec4} constitue la partie centrale de l'article, o\`u nous
allons prouver la r\'esurgence (locale) des solutions BKW formelles
\'el\'ementaires. La preuve se fait en deux \'etapes~:
\begin{enumerate}
\item La premi\`ere \'etape consiste \`a construire dans le cas o\`u la
  fonction $F$ est holomorphe au voisinage de l'origine
  (respectivement enti\`ere) des fonctions
  confluentes (respectivement fonctions confluentes r\'esurgentes)
  solutions de (\ref{eq4intro}) \`a support singulier la courbe
  alg\'ebrique $\mathcal{C} = \{(z, \xi),\,\, 9\xi^2 = 4z^3 \}$.
Cette construction repose essentiellement sur deux ingr\'edients : une
quantification de la transformation canonique associ\'ee \`a
l'op\'erateur principal intervenant dans l'\'equation
(\ref{eq4intro}), puis la r\'esolution d'une EDP singuli\`ere.\\
\item La deuxi\`eme \'etape consiste alors \`a d\'emontrer l'existence
  d'une famille de solutions BKW \'el\'ementaires qui peuvent
  \^etre vues comme la d\'ecomposition locale (respectivement
  d\'ecomposition) dans des germes de secteurs
  de Stokes (respectivement secteurs de Stokes) convenables des
  fonctions confluentes (respectivement fonctions confluentes
  r\'esurgentes) pr\'ec\'edemment construites.
\end{enumerate}

La section \ref{sec5} est consacr\'ee aux applications des r\'esultats
obtenus en section \ref{sec4}. Un premier paragraphe \'etablit
l'existence d'un th\'eor\`eme local r\'esurgent de r\'eduction tandis
qu'un deuxi\`eme paragraphe est consacr\'e \`a l'analyse BKW de
l'\'equation de Schr\"odinger (\ref{eq1intro}) induite par celle de
notre \'equation principale (\ref{eq4intro}). Un dernier paragraphe
expose quelques extensions possibles de nos r\'esultats.

Enfin, la section \ref{sec6} expose quelques pistes de recherche
d\'ecoulant naturellement de notre analyse.

Nous terminons par un appendice qui expose bri\`evement quelques notions
fondamentales utilis\'ees dans ce papier.

\subsection{Convention}

Dans l'analyse BKW exacte, tous les principaux objets
((pr\'e)sommation de Borel, secteurs de Stokes, etc\ldots) sont
relatifs \`a une direction donn\'ee $\alpha$, qui peut \^etre vue
comme un argument.\\
Dans tout ce qui va suivre, sauf mention contraire, nous supposerons
que $\alpha=0$, de sorte que $\Re (\varepsilon)>0 $
(et $|\varepsilon|$ assez petit).

\section{Cas de l'\'equation d'Airy}\label{sec2}

Nous nous concentrons ici sur l'\'equation d'Airy~:
\begin{equation}\label{eq3}
\frac{d^2 y}{d s^2} =\frac{s}{\varepsilon^2} y,
\end{equation}
c'est-\`a-dire sur l'\'equation (\ref{eq4intro}) lorsque $F=0$.\\
Comme nous l'avons dit, cette \'equation va nous servir de
r\'ef\'erence pour l'analyse BKW de l'\'equation (\ref{eq4intro}), du
fait que l'op\'erateur principal intervenant dans (\ref{eq4intro}) est
pr\'ecis\'ement celui d'Airy.\\
Nous rappelons ici les principaux r\'esultats connus concernant
l'analyse BKW de l'\'equation d'Airy.

\subsection{Aspect formel : le symbole BKW d'Airy}

Nous commen\c cons par introduire une solution BKW formelle "bien
normalis\'ee" associ\'ee \`a l'\'equation d'Airy~:
\begin{defn}
La solution BKW \'el\'ementaire suivante~:
\begin{equation}\label{wkb10}
 \left\{
 \begin{array}{l}
\displaystyle  A_{bkw}(z,\varepsilon) = 
\frac{e ^{\displaystyle  -\frac{2}{3}\frac{z^{3/2}}{\varepsilon} }}{z^{\frac{1}{4}}} \left(1 +
\sum_{n=1}^{+\infty} \alpha_n(z)\varepsilon^{n}\right) \\
\\
\displaystyle \alpha_n(z) = \left(-\frac{3}{4}
\right)^n\frac{\Gamma(n+\frac{1}{6})\Gamma(n+\frac{5}{6})}{2\pi \Gamma(n+1)} z^{-\frac{3n}{2}}  \hspace{10mm} n \geq 1. 
\end{array}
 \right.
\end{equation}
sera appel\'ee le {\em symbole BKW d'Airy}.
\end{defn}

Le symbole BKW d'Airy satisfait les propri\'et\'es fondamentales de r\'esurgence et
de sommabilit\'e (de Borel) suivantes~:

\begin{prop}\label{ressymbAiry}
Le  symbole BKW d'Airy est r\'esurgent
 sommable de  Borel en $\varepsilon^{-1}$, \`a d\'ependance r\'eguli\`ere en $z \neq 0$. 
\end{prop}

\subsection{Etude du ph\'enom\`ene de Stokes associ\'e}

\begin{figure}[thp]
\begin{center}
\begin{tabular}{ccc}
\includegraphics[width=2.1in]{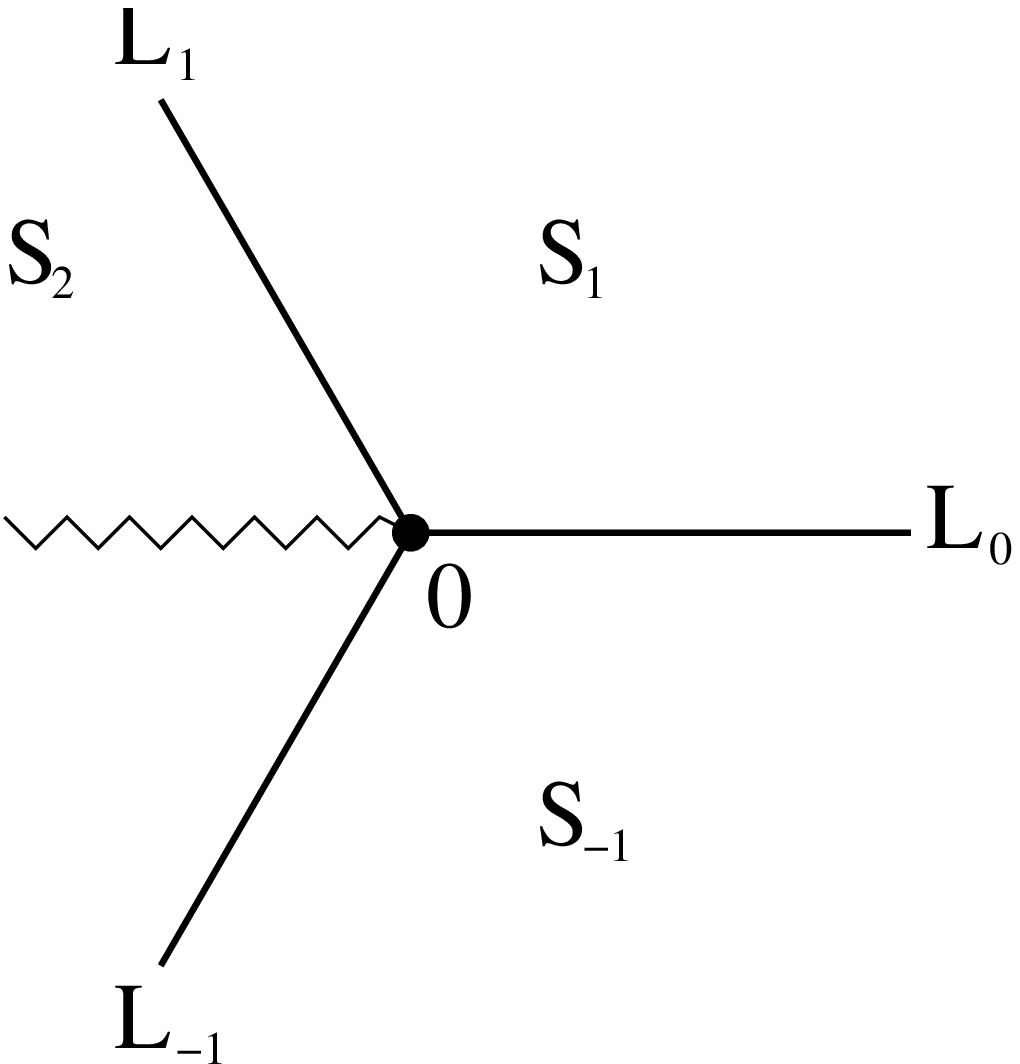} & \hspace{10mm} &
\includegraphics[width=1.5in]{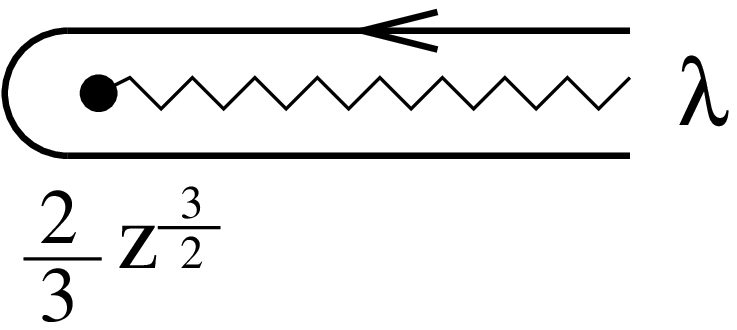} \\
Fig. \ref{fig:Stokes5}.a  &  &  Fig. \ref{fig:Stokes5}.b \\
\end{tabular}
\caption{ Fig. \ref{fig:Stokes5}.a  : $L_0$, $L_1$ et $L{-1}$ sont
les lignes de Stokes  
  (dans le $z$-plan) associ\'ees \`a la direction $\alpha=0$. 
Les trois secteurs de Stokes sont les secteurs ouverts connexes born\'es par les lignes de Stokes (en oubliant la ligne ondul\'ee). 
Fig. \ref{fig:Stokes5}.b : Le contour d'int\'egration dans le
$\xi$-plan. Les lignes ondul\'ees sont des coupures.
\label{fig:Stokes5}}
\end{center}
\end{figure}

Pour cette \'etude, nous renvoyons \`a \cite{Jidoumou, DDP93, DDP97,
  DP99} pour plus de d\'etails.\\
Rappelons ici que nous avons fait le choix de prendre la direction
$\alpha=0$ comme direction de sommation de Borel.\\
Les lignes de Stokes et les secteurs de Stokes
sont alors ceux dessin\'es sur la figure \ref{fig:Stokes5}.a.

Tant que $z$ reste dans l'un des secteurs de Stokes, le symbole BKW
d'Airy est sommable de Borel. 
Par exemple, fixons les conventions suivantes~:\\

\noindent {\bf Convention}: en dessinant une coupure comme sur la 
Fig. \ref{fig:Stokes5}.a, nous fixons la 
d\'etermination de $z^{3/2}$ ({\em resp.} $z^{1/4}$) de sorte que $z^{3/2}$ ({\em resp.}
$z^{1/4}$) est r\'eel positif le long de $L_0$. 
Nous notons $A_{bkw}^+(z, \varepsilon)$ la d\'etermination de $A_{bkw}(z, \varepsilon)$
ainsi d\'efinie, et $A_{bkw}^-(z, \varepsilon) := A_{bkw}^+(z, -\varepsilon)$.\\

\noindent {\bf Notation}:  nous avons vu dans la proposition
\ref{ressymbAiry} que le symbole BKW d'Airy $A_{bkw}^+$ est sommable
de Borel. \\
Nous noterons par~:
\begin{equation}\label{wkb11bis}
\mathcal{A}(z, \varepsilon)  =\mbox{\sc s}_0 \left(  A_{bkw}^+ \right) (z,\varepsilon) 
\end{equation} sa somme de Borel.\\  
Rappelons que cette derni\`ere est holomorphe en $(z, \varepsilon)$, $\Re (\varepsilon) >0$ et $z \in
S_1$ ({\em resp.} $S_{-1}$) et s'\'etend analytiquement en une fonction enti\`ere en $z$.
En particulier, $\displaystyle \mathcal{A}(z, \varepsilon)  = 2\sqrt{\pi}\varepsilon^{-1/6} Airy
(z\varepsilon^{-2/3})$, 
o\`u $Airy$ est la fonction d'Airy. \\

Historiquement, c'est par l'interm\'ediaire de l'\'equation d'Airy que
Stokes d\'ecouvrit le ph\'enom\`ene qui porte aujourd'hui son nom
(voir son article fondateur de 1857 \cite{Sto57}).\\ 
Il y a plusieurs fa\c cons de d\'ecrire le ph\'enom\`ene de Stokes :
le point de vue adopt\'e ici est de d\'ecrire ce ph\'enom\`ene comme une
rupture dans la d\'ecomposition de la fonction $ \mathcal{A}(z,
\varepsilon)$ lors de la travers\'ee d'une ligne de Stokes. 
Cette rupture est due \`a la pr\'esence de singularit\'es pour le
mineur associ\'e \`a  $\mathcal{A}(z, \varepsilon)$.\\
Pr\'ecisons les choses.
 
La sommabilit\'e de Borel induit une
correspondance bijective entre un d\'eveloppement formel et sa somme
de Borel de sorte que nous pouvons associer \`a $\mathcal{A}$ sa {\em d\'ecomposition}
$A_{bkw}^+$ pour $z \in S_1$ ({\em resp.} $S_{-1}$):
\begin{equation}\label{eqa5}
\begin{array}{c}
\mathcal{A}(z, \varepsilon)  \hspace{3mm}  \stackrel{\displaystyle \sigma_{S_1}}{\displaystyle
  \longrightarrow} \hspace{3mm}  A_{bkw}^+(z,\varepsilon).\\
\\
\left(\mbox{{\em resp.} } \mathcal{A}(z, \varepsilon)  \hspace{3mm}  
\stackrel{\displaystyle \sigma_{S_{-1}}}{\displaystyle
  \longrightarrow} \hspace{3mm}  A_{bkw}^+(z,\varepsilon). \right)
\end{array}
\end{equation}
 
Le fait que la d\'ecomposition de $\mathcal{A}(z, \varepsilon)$ dans
$S_1$ et $S_{-1}$ est donn\'ee par le m\^eme d\'eveloppement formel,
ou autrement dit, que la sommation de Borel et prolongement analytique en $z$
commutent encore lorsque l'on franchit la ligne de Stokes $L_0$, est d\^u au fait
que le symbole BKW d'Airy $A_{bkw}^+$ est {\em r\'ecessif  le long de $L_0$} (avec la
d\'etermination pr\'ec\'edemment choisie pour $z^{3/2}$). \\
En revanche, ce n'est plus vrai lorsque, venant de $S_1$  ({\em resp.} $S_{-1}$) l'on
traverse la ligne de Stokes $L_1$
({\em resp.} $L_{-1}$): pour $z$ sur ces lignes, un  ph\'enom\`ene de Stokes appara\^it,
et ce dernier est compl\`etement d\'ecrit par l'action de la
d\'erivation \'etrang\`ere suivante~:
\begin{equation}\label{wkb12}
\displaystyle  \dot \Delta_{-\frac{4}{3}z^{3/2}} A_{bkw}^+(z,\varepsilon) =
\ell A_{bkw}^+(z,\varepsilon) =-i A_{bkw}^-(z,\varepsilon)
\end{equation} 
o\`u $\ell$ est le prolongement analytique en $z$ autour de $0$
dans le sens trigonom\'etrique. Cela signifie que la d\'ecomposition de $\mathcal{A}$ pour
$z \in S_2$ (disons) devient:
\begin{equation}\label{eqa6}
\begin{array}{ccc}
\mathcal{A} (z, \varepsilon) & \stackrel{\displaystyle \sigma_{S_2}}{\displaystyle
  \longrightarrow} & A_{bkw}^+(z,\varepsilon) - \ell A_{bkw}^+(z,\varepsilon) = 
A_{bkw}^+(z,\varepsilon) +i A_{bkw}^-(z,\varepsilon) 
\end{array}
\end{equation} 
De m\^eme, pour $z \in L_0$, nous avons:
\begin{equation}\label{wkb12bis}
\displaystyle  \dot \Delta_{+\frac{4}{3}z^{3/2}} A_{bkw}^-(z,\varepsilon) =
\ell A_{bkw}^-(z,\varepsilon) = -i A_{bkw}^+(z,\varepsilon).
\end{equation} 

\begin{figure}[thp]
\begin{center}
\includegraphics[width=1.5in]{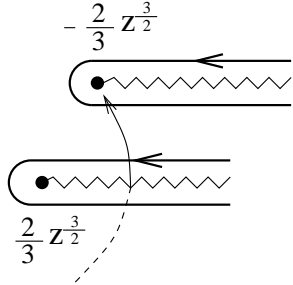} 
\caption{Effet du ph\'enom\`ene de Stokes
  d\'ecrit par (\ref{wkb12}) en termes de la d\'eformation du contour
  d'int\'egration pour la somme de Borel (\ref{wkb11}).
\label{fig:Stokes62}}
\end{center}
\end{figure}

La pr\'esence de ces deux singularit\'es (mobiles avec $z$) pour le
mineur associ\'e \`a $\mathcal{A}(z, \varepsilon)$ se traduit
\'egalement naturellement en termes de lieu singulier d'un majeur.\\ 
En effet, la somme de Borel de $A_{bkw}$ pour $z \in S_1$ (disons) peut \^etre d\'efinie comme une int\'egrale,
\begin{equation}\label{wkb11}
\mbox{\sc s}_0 \left(  A_{bkw} \right) (z,\varepsilon)  =  
\int_{\lambda}  e^{-\frac{1}{\varepsilon} \xi} \stackrel{\vee}{A_{bkw}}(z,\xi) \, d\xi.
\end{equation}
o\`u $\displaystyle \stackrel{\vee}{A_{bkw}}(z,\xi)$ est un majeur
associ\'e au symbole BKW d'Airy. Ce majeur est holomorphe sur le
rev\^etement universel de $\mathbb{C}^2 \backslash \mathcal{C}$, o\`u le
 support singulier $\mathcal{C}$ est la courbe alg\'ebrique $\mathcal{C} = \{(z, \xi),\,\, 9\xi^2 = 4z^3 \}$. Le
contour d'int\'egration $\lambda$ est dessin\'e sur la figure
\ref{fig:Stokes5}.b pour $z \in S_1$, et sa d\'eformation
pour $z \in S_2$  apr\`es la travers\'ee de la ligne de Stokes $L_1$
est dessin\'ee sur la figure \ref{fig:Stokes62}.

Notons pour terminer que la repr\'esentation int\'egrale (\ref{wkb11}) ci-dessus
peut \^etre d\'eduite de la repr\'esentation usuelle pour la fonction
d'Airy, plus pr\'ecis\'ement 
(\`a un facteur de normalisation pr\`es):
\begin{equation}\label{eq5bis}
\displaystyle  \int  e^{-\frac{1}{\varepsilon} S(z,\widehat{z})}  \, d\widehat{z} \hspace{5mm}
\mbox{o\`u} \hspace{5mm} S(z,\widehat{z}) =  z\widehat{z} -\frac{1}{3}\widehat{z}^3.
\end{equation}
Notre analyse dans la section \ref{sec4} sera bas\'ee sur une
extension de cette repr\'esentation int\'egrale.

\section{Analyse BKW formelle dans le cas g\'en\'eral}\label{sec3}

Nous nous focalisons maintenant sur l'\'equation~:
\begin{equation}\label{eq4}
\displaystyle  \frac{d^2 \Phi}{d z^2} - \frac{z}{\varepsilon^2} \Phi = F(z) \Phi,
\end{equation} 
en supposant d\'esormais que $F$ est une fonction analytique au voisinage
de l'origine quelconque.\\
Nous nous int\'eressons tout d'abord au probl\`eme de l'existence de
solutions BKW formelles de l'\'equation (\ref{eq4intro}) (de mani\`ere
analogue \`a la section \ref{sec2}).

\subsection{Existence de solutions BKW formelles}

Etant donn\'e que l'op\'erateur principal apparaissant dans
l'\'equation (\ref{eq4intro}) est celui d'Airy, il est naturel de
rechercher des solutions BKW formelles de la m\^eme forme que celle du
symbole BKW d'Airy.\\
Ceci nous conduit \`a la proposition suivante
(dont la d\'emonstration est imm\'ediate)~:
\begin{prop}
Il existe des solutions BKW formelles de l'\'equation (\ref{eq4intro})
de la forme~:
\begin{equation}\label{wkb1}
\Phi_{bkw}(z,\varepsilon) =  
\frac{e^{\displaystyle  -\frac{2}{3}\frac{ z^{3/2}}{\varepsilon} }}{z^{\frac{1}{4}}} (1 +
g_1(z) \varepsilon^{1}  + g_2(z) \varepsilon^{2} + \cdots).
\end{equation}
Dans ce cas, les fonctions $g_n$ v\'erifient les \'equations
(diff\'erentielles) de transport suivantes~:
\begin{equation}\label{wkb2}
 \left\{
 \begin{array}{l}
\displaystyle  32z^{5/2} \frac{dg_1}{dz} +16z^2F(z) -5 =0\\
\\
\displaystyle  32z^{5/2} \frac{dg_{n+1}}{dz} -16z^{2}\frac{d^2g_n}{dz^2}
+8z\frac{dg_n}{dz} + \left(16z^2F(z)-5\right)g_n  =0, \hspace{5mm} n
\geq 1.\\
\end{array}
 \right.
\end{equation}
\end{prop}

Bien \'evidemment, le d\'eveloppement (\ref{wkb1}), qui est
multivalu\'e en $z$, d\'epend du choix de la
d\'etermination pour  $z^{3/2}$ (de m\^eme que pour $z^{\frac{1}{4}}$).\\
Puisque l'\'equation (\ref{eq4}) est invariante sous l'action de $\varepsilon
\mapsto - \varepsilon$, nous en d\'eduisons que
\begin{equation}\label{wkb1bis}
\Phi_{bkw}(z,-\varepsilon)
\end{equation}
est une autre solution BKW formelle, et que de plus $\{ \Phi_{bkw}(z,\varepsilon),
\Phi_{bkw}(z,-\varepsilon) \}$ d\'efinit une base de solutions BKW formelles
pour l'\'equation  (\ref{eq4}). \\

\subsection{Solutions BKW \'el\'ementaires}

Nous voudrions obtenir une normalisation analogue \`a celle adopt\'ee
pour le symbole BKW d'Airy.\\
Pour cela, il est int\'eressant d'utiliser une autre
repr\'esentation de ces d\'eveloppements BKW. En \'ecrivant
$\Phi_{bkw}(z,\varepsilon)$ sous la forme 
\begin{equation}\label{wkb3}
\Phi_{bkw}(z,\varepsilon) = \exp \left(-\frac{1}{\varepsilon} \int^z P(t, \varepsilon) dt \right),
\end{equation}
l'\'equation (\ref{eq4}) devient:
\begin{equation}\label{wkb4}
\frac{1}{\varepsilon} \frac{dP}{dz} +\frac{1}{\varepsilon^{2}}\left(z-P^2\right) + F(z)=0.
\end{equation}
Cela signifie que si
\begin{equation}\label{wkb5}
P(z, \varepsilon) = \sum_{n \geq 0} p_n(z)\varepsilon ^{n}
\end{equation}
alors~:
\begin{equation}\label{wkb6}
 \left\{
 \begin{array}{l}
\displaystyle  p_0^2 =z\\
\displaystyle 2p_0p_1 = \frac{dp_0}{dz}\\
\displaystyle 2p_0p_2 = \frac{dp_1}{dz}-p_1^2+F(z)\\
\displaystyle 2p_0p_{n+1} = \frac{dp_{n}}{dz}-\sum_{1 \leq j \leq n}
p_jp_{n+1-j}, \hspace{10mm} n \geq 2. 
\end{array}
 \right.
\end{equation}
Nous montrons facilement par r\'ecurrence que~:
\begin{equation}\label{wkb6bis}
 \left\{
 \begin{array}{l}
\displaystyle  p_0(z) =z^{\frac{1}{2}}\\
\displaystyle p_1 (z)= \frac{1}{4z}\\
\\
\displaystyle p_n (z) \in z^{-\frac{3n-1}{2}}\mathbb{C}\{z\},  \hspace{10mm} n \geq 2. 
\end{array}
 \right.
\end{equation}
En introduisant la d\'ecomposition $\displaystyle  P = P_{pair} + P_{impair}$, 
$\displaystyle
 \left\{
 \begin{array}{l}
\displaystyle  P_{pair} = \sum_{k \geq 0} p_{2k}\varepsilon^{2k}\\
\displaystyle  P_{impair} = \sum_{k \geq 0} p_{2k+1}\varepsilon^{2k+1}
\end{array}
 \right.
$,
nous d\'eduisons de (\ref{wkb4}) que
$\displaystyle 
P_{impair} = \frac{\varepsilon}{2} \frac{P_{pair}^\prime}{P_{pair}}$ o\`u 
$\displaystyle P_{pair}^\prime = \frac{dP_{pair}}{dz}$.
Par cons\'equent, nous avons la repr\'esentation~:
\begin{equation}\label{wkb9}
\Phi_{bkw}(z,\varepsilon) =  
\frac{C(\varepsilon)}{\sqrt{P_{pair}(t, \varepsilon)}}\exp \left( -\frac{1}{\varepsilon}\int^z
  P_{pair}(t, \varepsilon) dt \right), \hspace{5mm} \mbox{avec } C(\varepsilon) \in
\mathbb{C}[[\varepsilon]].
\end{equation}

\begin{prop}\label{prop21}
Les solutions BKW formelles (\ref{wkb1}) 
de (\ref{eq4}) peuvent \^etre normalis\'ees de telle mani\`ere que
pour tout $n \geq 0$,
$\displaystyle g_n(z) \in z^{-\frac{3n}{2}}\mathbb{C}\{z\}$.
\end{prop}

\begin{proof}
Pour $n=1$, nous d\'eduisons de (\ref{wkb2}) que $g_1(z) = h_1(z) + Cste$,
o\`u $h_1(z) \in z^{-\frac{3}{2}}\mathbb{C}\{z\}$ tandis que $Cste$
est un nombre complexe quelconque. En choisissant $Cste=0$, cela
fournit le r\'esultat. \\
Maintenant, pour un $n \geq 1$ fix\'e, nous supposons que $\displaystyle g_n(z) \in
  z^{-\frac{3n}{2}}\mathbb{C}\{z\}$. De (\ref{wkb2}) nous tirons~:
$$ \displaystyle g_{n+1}(z) = -\frac{1}{32} \int H_n(z) \, dz
,$$ o\`u
$$H_n(z) = \frac{
  -16z^{2}g_n^{\prime \prime}(z)
+8zg_n^\prime (z) + \left(16z^2F(z)-5\right)g_n (z)}{z^{5/2}}
\in z^{-\frac{3n+5}{2}}\mathbb{C}\{z\}.$$
Si $n$ est pair, nous obtenons que $g_{n+1}(z) = h_{n+1}(z) + Cste$, o\`u 
$h_{n+1}(z) \in z^{-\frac{3n+3}{2}}\mathbb{C}\{z\}$. En choisissant
$Cste=0$ pour la constante d'int\'egration, cela donne le r\'esultat.
Si $n$ est impair, un $\ln(z)$ pourrait {\em a priori} appara\^itre par int\'egration,
mais cela serait en contradiction avec la repr\'esentation \'equivalente
(\ref{wkb9}) et la propri\'et\'e (\ref{wkb6bis}).
\end{proof}
 
\begin{defn}\label{forwkbsol}
Les solutions BKW formelles d\'ecrites dans la proposition \ref{prop21} seront
appel\'ees les {\em solutions BKW \'el\'ementaires de l'\'equation (\ref{eq4})}.
\end{defn}

\section{R\'esurgence des solutions BKW \'el\'ementaires}\label{sec4}

\subsection{Construction de fonctions confluentes}

Nous revenons maintenant aux solutions BKW \'el\'ementaires
d\'ecrites dans la proposition \ref{prop21}. Nous voudrions
``r\'ealiser'' le th\'eor\`eme de Borel-Ritt, c'est-\`a-dire
construire des fonctions analytiques dont l'asymptotique est
gouvern\'ee par (au moins une famille de) ces symboles BKW
\'el\'ementaires.

Le point de vue est donc ici "inverse" par rapport au cas d'Airy :
nous ne partons pas d'objets formels pour en d\'eduire des fonctions
analytiques par (pr\'e)sommation mais au contraire nous voulons
partir de fonctions confluentes (respectivement confluentes r\'esurgentes) et d\'eduire nos
objets formels (plus pr\'ecis\'ement une famille de symboles BKW \'el\'ementaires) par
d\'ecomposition dans des germes de secteurs de Stokes (respectivement
secteurs de Stokes) convenables. Ce point de vue "inverse" est en
effet souvent plus commode lorsque l'on manipule des objets
d\'ependant analytiquement d'un param\`etre (typiquement lorsqu'on
\'etudie la r\'esurgence param\'etrique d'objets formels).

\subsubsection{Repr\'esentation de type Laplace}

Puisque le symbole principal $p^2-z$ de l'op\'erateur d\'efinissant
l'\'equation (\ref{eq4}) est simplement l'op\'erateur d'Airy, en nous
inspirant des deux diff\'erentes repr\'esentations de la somme de
Borel du symbole BKW d'Airy, nous pouvons rechercher de telles
solutions analytiques sous deux formes~:
\begin{enumerate}
\item Une premi\`ere piste est de partir de la repr\'esentation (\ref{eq5bis}) ci-dessus, en
pensant $\displaystyle S(z,\widehat{z}) =  z\widehat{z} -\frac{1}{3}\widehat{z}^3$ comme une
fonction g\'en\'eratrice de la transformation canonique $(p,z)
\leftrightarrow (\widehat{p}, \widehat{z})$ dans l'espace cotangent, dont
l'effet est de redresser la sous-vari\'et\'e Lagrangienne $\widehat{p}
=p^2-z=0$.\\
Cette piste de recherche nous am\`ene \`a consid\'erer, comme dans
\cite{Ph88}, la quantification de la transformation canonique, i.e
rechercher des solutions de la forme~:
\begin{equation}\label{eq6bis}
 \displaystyle  \Phi(z,\varepsilon) =  \int  e^{-\frac{1}{\varepsilon} S(z,\widehat{z})} {\widetilde \varphi}(\widehat{z},
\varepsilon) \, d\widehat{z}.
\end{equation}
\item Une seconde piste est de rechercher des solutions de (\ref{eq4})
d\'efinies comme somme de Borel, i.e~:
\begin{equation}\label{eq6}
\displaystyle \Phi(z,\varepsilon) =  \int  e^{-\frac{1}{\varepsilon} \xi} \stackrel{\vee}{\Phi}(z,\xi) \, d\xi,
\end{equation}
o\`u $\displaystyle \stackrel{\vee}{\Phi}(z,\xi)$ doit \^etre un
majeur d'une microfonction confluente convenable (au sens
d\'evelopp\'e dans l'appendice \ref{sec7}). Ce que nous entendons
par ``convenable'' est la chose suivante : dans la repr\'esentation
int\'egrale (\ref{eq6}), en d\'erivant sous le signe somme et en
int\'egrant formellement par parties, nous traduisons le fait que 
$\Phi$ est solution de (\ref{eq4}) par le fait  de demander \`a 
 $\stackrel{\vee}{\Phi}$ de satisfaire l'\'equation: 
\begin{equation}\label{eq8}
\frac{\partial^2 \stackrel{\vee}{\Phi}}{\partial z^2} - z \frac{\partial^2
  \stackrel{\vee}{\Phi}}{\partial \xi^2} = F(z) \stackrel{\vee}{\Phi}.
\end{equation}
\end{enumerate}

Au lieu de rechercher directement des solutions pour l'EDP
(\ref{eq8}), nous allons combiner les deux id\'ees pr\'ec\'edentes
li\'ees aux repr\'esentations int\'egrales (\ref{eq6bis}) et
(\ref{eq6}).\\
En faisant dans (\ref{eq6}) le changement de variable $\xi
\leftrightarrow \widehat{z}$  d\'efini par $\xi = S(z, \widehat{z})$, nous
obtenons la repr\'esentation int\'egrale~:
\begin{equation}\label{eq7}
\displaystyle \Phi(z,\varepsilon) =  
\int_\gamma  e^{-\frac{1}{\varepsilon} S(z, \widehat{z})} \Psi(z,\widehat{z})  \, d\widehat{z}, \hspace{10mm}
\stackrel{\vee}{\Phi}(z,\xi)\big|_{\displaystyle \xi = S(z, \widehat{z})} = \frac{\Psi(z,\widehat{z})}{z-\widehat{z}^2},
\end{equation}
o\`u le chemin d'int\'egration $\gamma$ est, pour l'instant, vu comme
un chemin sans fin, allant \`a l'infini dans les zones o\`u $\displaystyle \Re \big( \frac{1}{\varepsilon}
S(z,\widehat{z}) \big) \rightarrow +\infty$.\\
En posant $\tpsi(z,\z) = \stackrel{\vee}{\Phi} (z, \xi)$, nous
d\'eduisons facilement de (\ref{eq8}) que $\tpsi$ doit \^etre solution
de l'EDP lin\'eaire suivante~:
\begin{equation}\label{eq9}
\frac{\partial^2 \tpsi}{\partial z^2} - \frac{2\z}{z-\z^2}\frac{\partial^2
  \tpsi}{\partial z \partial \z} - \frac{1}{z-\z^2}\frac{\partial^2 \tpsi
 }{\partial \z^2} = F(z) \tpsi.
\end{equation}

\subsubsection{R\'esolution de l'EDP singuli\`ere associ\'ee}\label{sec412}

\paragraph{Deux exemples}

\begin{itemize}
\item Lorsque $F(z) = \lambda^2$, $\lambda \in \mathbb{C}$, nous
  montrons facilement que 
\begin{equation}
\tpsi (z, \z) = \frac{e^{\pm \lambda (z-\z^2)}}{z-\z^2}
\end{equation}
sont solutions particuli\`eres de (\ref{eq9}). En choisissant dans
(\ref{eq7}) des contours d'int\'egration convenables $\gamma$, nous
obtenons ainsi une base de solutions r\'esurgentes de (\ref{eq4}) de
la forme~:
\begin{equation}
\Phi(z,\varepsilon) = \int_\gamma  e^{-\frac{1}{\varepsilon} S(z, \z)} \Psi(z,\z)
\, d\z, \hspace{10mm} \Psi(z,\z) = e^{\pm \lambda (z-\z^2)}.
\end{equation}

\item Lorsque 
$F(z) = \lambda^2 z$, $\lambda \in \mathbb{C}$, des solutions particuli\`eres
de (\ref{eq9}) sont donn\'ees par
\begin{equation}
\displaystyle \tpsi (z, \z) = \frac{e^{\pm \frac{1}{3}\lambda (z-\z^2)\sqrt{4z-\z^2}}}{z-\z^2},
\end{equation}
qui, par lin\'earit\'e, fournissent des solutions de (\ref{eq4}) de la
forme~:
\begin{equation}
\Phi(z,\varepsilon) = \int_\gamma  e^{-\frac{1}{\varepsilon} S(z, \z)} \Psi(z,\z)  \, d\z,
\hspace{10mm} 
\Psi(z,\z) = \cosh \left( 
\frac{1}{3}\lambda (z-\z^2)\sqrt{4z-\z^2}\right).
\end{equation}
Notons que $\Psi(z,\z)$ est holomorphe pour
$(z,\z) \in \mathbb{C}^2$, et que cette int\'egrale converge pour $|\varepsilon
|$ assez petit (pour un choix convenable de $\gamma$).
\end{itemize}

\paragraph{R\'esolution dans le cas g\'en\'eral}

Dans la repr\'esentation int\'egrale (\ref{eq7}), ayant en t\^ete la
m\'ethode du col, nous demandons \`a la fonction $\Psi(z, \z)$ d'\^etre
holomorphe au voisinage du lieu $\displaystyle  \frac{\partial
    S(z,\z)}{\partial \z} = 0$ d\'efinissant les points cols.
Puisque
$\displaystyle  \frac{\partial
    S(z,\z)}{\partial \z} = z-\z^2$, nous introduisons la transformation~:
\begin{equation}\label{eq18}
 \left\{
 \begin{array}{l}
\displaystyle (z, \z) \leftrightarrow (z,x= z-\z^2) \\
\\
\displaystyle 
\psi(z,x) := \Psi(z,\z)  = \tpsi(z,\z) (z-\z^2).     
\end{array}
 \right.
\end{equation}  
Par cette transformation, l'\'equation (\ref{eq9}) se traduit pour
$\psi$ en l'\'equation suivante~:
\begin{equation}\label{eq19}
x^2\frac{\partial^2 \psi}{\partial x^2} + (4xz-2x^2)\frac{\partial^2
  \psi}{\partial x \partial z} +x^2\frac{\partial^2
  \psi}{\partial z^2} +(2x-4z)\frac{\partial \psi}{\partial z} -x^2
F(z) \psi =0.
\end{equation}
Nous allons maintenant rechercher des solutions holomorphes de
l'\'equation (\ref{eq19}) pour $x$ au voisinage de
z\'ero (et $z$ proche de $0$ \'egalement).\\
Etant donn\'e que dans l'\'equation (\ref{eq19}), $x=0$ est un point
singulier, le r\'esultat est non trivial car il ne peut d\'ecouler
simplement du th\'eor\`eme de Cauchy-Kovalevska.\\

Nous allons d'abord commencer par regarder l'existence de solutions
formelles pour (\ref{eq19}) de la forme 
\begin{equation}\label{eq20}
\psi(z,x) = \sum_{n \geq 0} a_n(z) x^n.
\end{equation}
 
\begin{lem}\label{lem2}
Soit $h(z)$ une fonction holomorphe au voisinage de l'origine. Alors 
il existe un unique d\'eveloppement formel $\displaystyle \psi(z,x)
= \sum_{n \geq 0} a_n(z) x^n$ solution de (\ref{eq19}) tel que les
$a_n (z)$ soient des fonctions holomorphes au voisinage de $z=0$, avec
\begin{equation}\label{eq24}
 \left\{
\begin{array}{l}
a_0 (z)=1  \\
a_1(z)=h(z).
\end{array}
 \right. 
\end{equation}
Dans ce cas, nous avons de plus, pour $n \geq 2$~: 
\begin{equation}\label{eq25}
a_n(z) = \frac{1}{n-1} \int_0^1 u^{n-1}\bigg(-a_{n-2}^{\prime
    \prime}(u^4z)+2(n-2)a_{n-1}^{\prime}(u^4z)+F(u^4z)a_{n-2}(u^4z)\bigg) \, du
\end{equation} 
\end{lem}

\begin{proof}
En rempla\c cant $\psi(z,x)$ par (\ref{eq20}) dans l'\'equation
(\ref{eq19}), et en identifiant les puissances de $x$, nous obtenons
le syst\`eme suivant~:
\begin{equation}\label{eq23}
 \left\{
 \begin{array}{l}
\displaystyle \frac{\partial a_0 }{\partial z}=0  \\
\\
\displaystyle 
4z\frac{\partial a_n }{\partial z}+na_n = \frac{1}{n-1}\left(
-\frac{\partial^2 a_{n-2} }{\partial z^2} + 2(n-2)\frac{\partial a_{n-1} }{\partial z}+F(z)a_{n-2}\right), \, \, \,  \mbox{ pour } n \geq 2.
\end{array}
 \right.
\end{equation} 
Il suffit alors d'int\'egrer l'\'equation (\ref{eq23}) en tenant compte
du fait que les $a_n$ doivent \^etre holomorphes au voisinage de
$z=0$.
\end{proof}

Nous avons d\'emontr\'e au lemme \ref{lem2} l'existence d'une famille de
solutions formelles de (\ref{eq19}). \\
A notre connaissance, les th\'eories classiques (voir \cite{GT96, KKS79}) pour analyser la
convergence de ces solutions formelles de l'EDP singuli\`ere
(\ref{eq19}) ne s'appliquent pas dans notre cas. \\
Afin de montrer la convergence, nous allons utiliser le r\'esultat suivant~:

\begin{lem}\label{intrepre}
La s\'erie formelle donn\'ee au lemme \ref{lem2}
repr\'esente la s\'erie de Taylor d'une fonction holomorphe
  $\psi(z,x)$ au voisinage de $(z,x)=(0,0)$ si et seulement si 
$\displaystyle \varphi(z,x) := \frac{\psi(z,x)}{x} -
\frac{1}{x}$  satisfait l'\'equation int\'egrale suivante~:
\begin{equation}\label{le0ter}
\begin{array}{l}
\displaystyle \varphi(z,x)  =  h(z) +  \int_0^1 du \int_0^{ux} dt
F(zu^4)+2x\int_0^1 u\partial_1\varphi(zu^4,ux) du
 \\
\displaystyle \hspace{20mm} -  \int_0^1 du \int_0^{ux} dt   \Big( 
 t\partial_1^2 \varphi(zu^4,t) +2\partial_1\varphi(zu^4,t)-t F(zu^4)
 \varphi(zu^4,t)\Big),
\end{array}
\end{equation}
o\`u $\displaystyle \partial_1^2 \varphi 
:= \frac{\partial^2 \varphi}{\partial z^2}$ et $\displaystyle \partial_1 \varphi 
:= \frac{\partial \varphi}{\partial z}$.
\end{lem} 

\begin{proof}
Nous consid\'erons la s\'erie $\displaystyle \psi(z,x) = \sum_{n
  \geq 0} a_n(z) x^n$ donn\'ee par le lemme \ref{lem2}, en supposant la convergence.\\
Comme  $\displaystyle \int_0^x x_1^{n-2} \, dx_1 =
\frac{x^{n-1}}{n-1}$, nous pouvons \'ecrire, pour 
 $u \in [0,1]$ et $(z,x)$ dans un voisinage de l'origine,
$$ \begin{array}{ll}
\displaystyle -\sum_{n \geq 2} \frac{u^{n-1}}{n-1}a_{n-2}^{\prime \prime}(zu^4) x^n
 & = 
\displaystyle -x \sum_{n \geq 2} u^{n-1} a_{n-2}^{\prime \prime}(zu^4)
\int_0^x x_1^{n-2} \, dx_1 \\
\\
& \displaystyle  = -xu \int_0^x \sum_{n \geq 2} a_{n-2}^{\prime \prime}(zu^4)
(ux_1)^{n-2} \, dx_1 \\
\\
& \displaystyle  = -xu \int_0^x \sum_{n \geq 0} a_{n}^{\prime \prime}(zu^4)
(ux_1)^{n} \, dx_1 \\
\\
& \displaystyle  = -xu \int_0^x  \partial_1^2 \psi(zu^4,ux_1) \, dx_1 \\
\\
& \displaystyle  = -x \int_0^{ux}  \partial_1^2 \psi(zu^4,t) \, dt,
\end{array}
$$
o\`u nous avons utilis\'e $a_0(z) = 1$ (voir (\ref{eq24})). 
Par cons\'equent,
\begin{equation}\label{le1}
\begin{array}{ll}
\displaystyle \sum_{n \geq 2} \left( \int_0^1 \left( 
\frac{u^{n-1}}{n-1}-a_{n-2}^{\prime \prime}(zu^4) \right) du \right) x^n
 & 
\displaystyle  = -x \int_0^1 du \int_0^{ux}  \partial_1^2 \psi(zu^4,t) \, dt.
\end{array}
\end{equation}
Par ailleurs, nous avons~:
$$
\begin{array}{ll}
\displaystyle \sum_{n \geq 2} \frac{u^{n-1}}{n-1}F(zu^4)a_{n-2}(zu^4) x^n
 & = 
\displaystyle x \sum_{n \geq 2} u^{n-1}F(zu^4)a_{n-2}(zu^4)
\int_0^x x_1^{n-2} \, dx_1 \\
\\
& \displaystyle  = ux  \int_0^x \sum_{n \geq 2}F(zu^4) a_{n-2}(zu^4)(ux_1)^{n-2} \, dx_1 \\
\\
& \displaystyle  = ux  \int_0^x \sum_{n \geq 0}F(zu^4) a_{n}(zu^4)(ux_1)^{n} \, dx_1 \\
\\
& \displaystyle  = ux \int_0^x  F(zu^4) \psi(zu^4, ux_1) \, dx_1 \\
\\
& \displaystyle  = x \int_0^{ux}  F(zu^4) \psi(zu^4,t) \, dt, \\
\end{array}
$$
de sorte que
\begin{equation}\label{le2}
\begin{array}{ll}
\displaystyle \sum_{n \geq 2} \left(
 \int_0^1 \left(  \frac{u^{n-1}}{n-1}F(zu^4)a_{n-2}(zu^4)\right) du \right) x^n
 &  \displaystyle  =  x \int_0^1 du \int_0^{ux}  F(zu^4)\psi(zu^4,t) \, dt.
\end{array}
\end{equation}
Enfin, nous avons
$$ \begin{array}{ll}
\displaystyle \sum_{n \geq 2} 2(n-2)\frac{u^{n-1}}{n-1}a_{n-1}^{\prime}(zu^4) x^n
& = 
\displaystyle \sum_{n \geq 2} 2u^{n-1}a_{n-1}^{\prime}(zu^4) x^n -
2x\sum_{n \geq 2} a_{n-1}^{\prime}(zu^4)u^{n-1}\frac{x^{n-1}}{n-1} \\
\\
 & = 
\displaystyle 2x\sum_{n \geq 0}a_{n}^{\prime}(zu^4) (ux)^n - 2x\sum_{n \geq 2} u^{n-1} a_{n-1}^{\prime}(zu^4)
\int_0^x x_1^{n-2} \, dx_1 \\
\\
& \displaystyle  = 2x\partial_1\psi(zu ^4,ux)-2x \int_0^x \sum_{n \geq 2} a_{n-1}^{\prime}(zu^4)
(ux_1)^{n-1} \, \frac{dx_1}{x_1} \\
\end{array}
$$
et par suite
$$ \begin{array}{ll}
\displaystyle \sum_{n \geq 2} 2(n-2)\frac{u^{n-1}}{n-1}a_{n-1}^{\prime}(zu^4) x^n
& 
\displaystyle  = 2x\partial_1\psi(zu ^4,ux)-2x \int_0^x \sum_{n \geq 0} a_{n}^{\prime}(zu^4)
(ux_1)^{n} \, \frac{dx_1}{x_1} \\
\\
& \displaystyle  = 2x\partial_1\psi(zu ^4,ux)-2x \int_0^x
\partial_1\psi(zu^4,ux_1) \, \frac{dx_1}{x_1} \\
\\
& \displaystyle  = 2x\partial_1\psi(zu^4,ux)-2x\int_0^{ux}  \partial_1 \psi(zu^4,t) \, \frac{dt}{t},
\end{array}
$$
o\`u nous avons utilis\'e $a_0(z) = 1$ (voir (\ref{eq24})). 
Par cons\'equent, \\
$$\displaystyle \sum_{n \geq 2} \left( \int_0^1 \left( 
\frac{u^{n-1}}{n-1}2(n-2)a_{n-1}^{\prime}(zu^4) \right) du \right) x^n$$
\begin{flushright}
\begin{equation}\label{le3}
= 2x\int_{0}^{1} \partial_1\psi(zu^4,ux) du - 2x \int_0^1  \int_0^{ux}  \partial_1 \psi(zu^4,t) \, \frac{dt}{t}du.
\end{equation}
\end{flushright}
Maintenant en utilisant (\ref{eq24}) et (\ref{eq25}), nous d\'eduisons
de (\ref{le1}), (\ref{le2}) et (\ref{le3}) que~: 
$$
\begin{array}{l}
\displaystyle \qquad \qquad  \psi(z,x)   = 1+h(z)x + \sum_{n \geq 2} a_n(z) x^n \\
\displaystyle  =   1+h(z)x-x \int_0^1 du \int_0^{ux}  \partial_1^2
\psi(zu^4,t) \, dt+ x \int_0^1 du \int_0^{ux}  F(zu^4)\psi(zu^4,t) \,
dt
\end{array} $$
\begin{flushleft}
$$ \displaystyle \qquad \qquad +2x\int_{0}^{1} \partial_1\psi(zu^4,ux) du - 2x \int_0^1 du \int_0^{ux}  \partial_1 \psi(zu^4,t) \,\frac{dt}{t}.$$
\end{flushleft}
En se rappelant que $\displaystyle \varphi(z,x) := \frac{\psi(z,x)}{x} -
\frac{1}{x}$, cela nous donne (\ref{le0ter}).
\end{proof}

Le lemme \ref{intrepre} va nous permettre de prouver la convergence
des d\'eveloppements formels d\'efinis dans le lemme \ref{lem2}.
A cet effet, introduisons une d\'efinition.

\begin{defn}
Si $W$ est un ouvert born\'e de $\mathbb{C}^n$,
  $n \geq 1$,  et $E$ espace de Banach, nous notons par
$H(\overline{W},E)$ l'espace des fonctions 
$f : Z \mapsto f(Z) \in E$ qui sont continues pour $Z \in \overline{W}$
et holomorphes dans $W$. 
\end{defn}

Nous rappelons le r\'esultat classique suivant~:

\begin{prop}
Soit $W$ un ouvert born\'e de $\mathbb{C}^n$, $n \geq 1$, et 
$E$ un espace de Banach. Nous munissons l'espace $H(\overline{W},E)$
de la norme du maximum~:
$$\|f \|_W = \sup_{Z \in W} |f(Z)|.$$
Alors $(H(\overline{W},E), \|. \|_W)$ est un espace de Banach.
\end{prop}

Dans toute la suite, $D(0,l) \subset \mathbb{C}$ d\'esigne le disque
ouvert centr\'e en $0$ de rayon $l>0$.

\begin{thm}\label{prop22}
Soit $R>0$, $r_1>0$ et $r_0$ tel que $0<r_0<r_1$.
Notons $d_0=r_1-r_0$.\\
Supposons que $F \in H(\overline{D(0,r_1)},\mathbb{C})$ et 
$h \in H(\overline{D(0,r_1)},\mathbb{C})$.
Posons $\displaystyle r^\prime =\min \bigg\{\frac{3r_1}{2e}\bigg(-1+\sqrt{1+\frac{4r_0d_0}{9er_1^2}}\bigg),  R\bigg\}$.\\
Alors, il existe une unique fonction holomorphe 
$\psi(z,x)$ au voisinage de $(0,0)$ solution de l'\'equation (\ref{eq19}), et
satisfaisant les conditions initiales suivantes~:
\begin{equation}
 \left\{
\begin{array}{l}
\displaystyle \psi(z,0)=1  \\
\\
\displaystyle \frac{\partial \psi}{\partial x}(z,0) = h(z).
\end{array}
 \right. 
\end{equation}
De plus, $\psi(z,x)$ s'\'etend analytiquement sur $D(0,r_1) \times D(0,r^\prime)$.
\end{thm}

\begin{proof} 
La preuve est inspir\'ee plus ou moins de techniques standards (voir,
par exemple, \cite{Treves}, \S 17).
\begin{enumerate}
\item  Pour $0 \leq s \leq 1$ nous notons
$U_s := D(0,r_s)$ avec $r_s := r_0+sd_0$. 
\item Pour $u \in [0,1]$ et $x \in \overline{D(0,R)}$ nous introduisons les fonctions
  $\displaystyle T: (u,x) \mapsto T{(u,x)}$ et $\displaystyle L: (u,x) \mapsto L{(u,x)}$ d\'efinies par
\begin{equation}\label{pre2}
\left\{
\begin{array}{l}
\displaystyle  T{(u,x)} : \psi(z) \mapsto \partial^1 \psi(zu^4)\\
\\
\displaystyle  L{(u,x)} : \psi(z) \mapsto x\partial^2 \psi(zu^4) +2T(u,x)\psi(z)- x F(zu^4)\psi(zu^4).
\end{array}
\right.
\end{equation} 
Nous voyons $T(u,x)$ et $L(u,x)$ comme des op\'erateurs lin\'eaires
agissant sur l'espace de Banach
$H(\overline{U_s}, \mathbb{C})$ et \`a valeurs dans $H(\overline{U_{s^\prime}},\mathbb{C})$,
$\displaystyle T{(u,x)},L{(u,x)} :  H(\overline{U_s},\mathbb{C}) \rightarrow   H(\overline{U_{s^\prime}},\mathbb{C})$,
o\`u  $0\leq s^\prime < s \leq 1$.\\
Par les formules de Cauchy,
\begin{equation}
\left\{ 
\begin{array}{l}
\displaystyle \partial^1  \psi(zu^4) = \frac{1}{2 i \pi}
\oint  \frac{\psi(t)}{(t-zu^4)^{2}} dt \\ 
\\
\displaystyle \partial^2  \psi(zu^4) = \frac{2}{2 i \pi}
\oint  \frac{\psi(t)}{(t-zu^4)^{3}} dt,
\end{array} 
\right.
\end{equation}
o\`u nous int\'egrons dans le sens direct sur un cercle centr\'e en
$zu^4$. Etant donn\'e que $\displaystyle r_s-u^4 r_{s^\prime} = (1-u^4)r_0+  (s-u^4
s^\prime) d_0$, nous avons~:
\begin{equation}
\left\{ 
\begin{array}{l}
\displaystyle \|\partial^1  \psi_{(u)}\|_{u_{s'}} \leq
\frac{\|\psi\|_{u_s}}{(1-u^4)r_0+(s-u^4s')d_0} \\ 
\\
\displaystyle \|\partial^2  \psi_{(u)}\|_{u_{s'}} \leq
\frac{2\|\psi\|_{u_s}}{\big((1-u^4)r_0+(s-u^4s')d_0\big)^2},
\end{array} 
\right.
\end{equation}
o\`u $\displaystyle \partial^i  \psi_{(u)} : z \mapsto   
\partial^i  \psi (zu^4)$, avec $i=1,2$.
Par suite, pour tout $(u,x) \in [0,1] \times \overline{D(0,R)}$, nous avons~:
\begin{equation}\label{pre3}
\left\{ 
\begin{array}{l}
\displaystyle \|T(u,x)\psi\|_{u_{s'}} \leq
\frac{\|\psi\|_{u_s}}{(1-u^4)r_0+(s-u^4s')d_0} \\ 
\\
\displaystyle \|L(u,x)\psi\|_{u_{s'}} \leq
\frac{2|x|\|\psi\|_{u_s}}{\big((1-u^4)r_0+(s-u^4s')d_0\big)^2}+\frac{2\|\psi\|_{u_s}}{(1-u^4)r_0+(s-u^4s')d_0}\\
\\ \qquad \qquad \qquad \qquad + |x|\|F\|_{D(0,r_1)}\|\psi\|_{u_s} 
\end{array} 
\right.
\end{equation}

\item Introduisons maintenant~:
\begin{equation}\label{pre4}
\left\{
\begin{array}{l}
\displaystyle \theta_0(x):= z \mapsto  h(z) +  \int_0^1 du \int_0^{ux}
F(zu^4) \, dt\\
\\
\displaystyle \theta_{k+1}(x) =  \theta_0(x) + 2x\int_0^1
uT(u,x)\theta_k(ux) \,du-\int_0^1 du
\int_0^{ux} L{(u,t)} \theta_k (t) \, dt, \hspace{1mm} k \geq 0.
\end{array}
\right. 
\end{equation}
Evidemment (\ref{pre4}) d\'efinit une suite 
$\left( \theta_k \right)_k$ de fonctions holomorphes en $x \in
D(0,R)$, continues pour $x \in \overline{D(0,R)}$, \`a valeurs dans 
$H(\overline{U_s},\mathbb{C})$,i.e pour tout $0 \leq s < 1$:
\begin{equation}\label{pre55}
\forall k \geq 0, \, \theta_k \in H \left( \overline{D(0,R)}, H(\overline{U_s},\mathbb{C}) \right).
\end{equation}
Posons \'egalement~:
\begin{equation}\label{pre5}
\left\{ 
\begin{array}{l}
\displaystyle \delta_0(x) :=  \theta_0(x)\\
\\
\displaystyle \delta_{k+1}(x) :=  \theta_{k+1}(x) -
\theta_{k}(x) \\
\\
\displaystyle \qquad \quad \, \, \,= 2x\int_0^1 uT(u,x)\delta_k(ux) \, du-\int_0^1 du \int_0^{ux} L{(u,t)} \delta_k (t) \, dt , k \geq 0.
\end{array}
\right.
\end{equation}
Observons dans un premier temps que, pour tout $0 \leq s <  1$, et tout $x  \in \overline{D(0,R)}$,
\begin{equation}\label{pre6}
\| \delta_0 (x)  \|_{U_{s}} \leq M, \hspace{5mm}\mbox{ avec } M = \|h\|_{D(0,r_1)} + \frac{R}{2}\|F\|_{D(0,r_1)}.
\end{equation}
Nous allons alors montrer le lemme suivant
\begin{lem}\label{lemfon}
Pour tout $\displaystyle 0 \leq  s < 1$, pour tout $k \in \mathbb{N}$
et tout $x \in \overline{D(0,R)}$,
\begin{equation}\label{pre7}
\| \delta_k (x)  \|_{U_{s}} \leq M \Big(\frac{e|x|\big(\alpha_k|x|+\beta\big)}{r_0d_0(1-s)}\Big)^k,
\end{equation}
avec $\displaystyle \alpha_k=1+\frac{r_0d_0(1-s)\|F\|_{D(0,r_1)}}{k}$ et $\beta=3r_1$.
\end{lem}

\begin{proof}
Nous proc\'edons par r\'ecurrence sur $k$.\\ 
Le cas $k=0$ est donn\'e par (\ref{pre6}).\\
Supposons maintenant que (\ref{pre7}) soit satisfaite pour un $k \in
\mathbb{N}$ donn\'e et pour tout $\displaystyle
0 \leq  s < 1$. \\
Pour tout  $\displaystyle 0 \leq s^\prime < s < 1$, nous d\'eduisons
de (\ref{pre5}) et (\ref{pre3}) que~:
$$ \| \delta_{k+1} (x) \|_{U_{s^\prime}} \leq 
2|x|\int_0^1\frac{u\|\delta_{k} (ux)
  \|_{U_{s}}}{(1-u^4)r_0+(s-u^4s^{\prime})d_0} \, du $$
\begin{flushright}
$$ + \int_0^1 du \int_0^{u|x|} 
\Bigg( \frac{2t}{\big(
(1-u^4)r_0+  (s-u^4s^\prime) d_0\big)^2} + \frac{2}{
(1-u^4)r_0+  (s-u^4s^\prime) d_0} $$
\end{flushright} 
\begin{flushright}
\begin{equation}\label{pre100}
 + t\|F\|_{D(0,r_1)}\Bigg) 
\|\delta_{k} (t) \|_{U_{s}} \, dt.
\end{equation}
\end{flushright}
Par l'hypoth\`ese de r\'ecurrence faite sur $\delta_k(x)$, nous avons alors~:

$$ \begin{array}{lll}
I_1 & = &\displaystyle 2|x|\int_0^1\frac{u\|\delta_{k} (ux)
  \|_{U_{s}}}{(1-u^4)r_0+(s-u^4s^{\prime})d_0} \, du \\
\\   
& \leq & \displaystyle
\frac{2Me
  ^k}{\big(r_0d_0(1-s)\big)^k}\int_0^1\frac{(u|x|)^{k+1}(\alpha_ku|x|+\beta)^k}{(1-u^4)r_0+(s-u^4s^{\prime})d_0} \, du\\
\\
&  \leq &\displaystyle \frac{2Me
  ^k|x|^{k+1}(\alpha_k|x|+\beta)^k}{\big(r_0d_0(1-s)\big)^k}\int_0^1\frac{u^{k+1}}{(1-u^4)r_0+(s-u^4s^{\prime})d_0} \, du.
\end{array} $$ 
Or, nous avons la majoration suivante~:
$$ \displaystyle
\int_0^1\frac{u^{k+1}}{(1-u^4)r_0+(s-u^4s^{\prime})d_0} du \leq
\frac{1}{d_0(s-s^{\prime})}\int_0^1u^{k+1} du \leq \frac{1}{(k+1)d_0(s-s^{\prime})}.$$
Par suite, nous en d\'eduisons que~:
\begin{equation}\label{IneqI1}
I_1 \leq
\frac{Me^k(\alpha_k|x|+\beta)^k|x|^{k+1}}{(k+1)\big(r_0d_0(1-s)\big)^k}\Big(\frac{2r_1}{r_0d_0(s-s^{\prime})}\Big).
\end{equation}

De m\^eme, nous avons~:
$$ \begin{array}{lll}
I_2 & = &\displaystyle \int_0^1 du \int_0^{u|x|}\frac{2t\|\delta_{k} (t)
  \|_{U_{s}}}{\big((1-u^4)r_0+(s-u^4s^{\prime})d_0\big)^2} \, dt \\
\\   
& \leq & \displaystyle
\frac{2Me
  ^k(\alpha_k|x|+\beta)^k}{\big(r_0d_0(1-s)\big)^k}\int_0^1 du
\int_0^{u|x|} \frac{t^{k+1}}{\big((1-u^4)r_0+(s-u^4s^{\prime})d_0\big)^2} \, dt\\
\\
&  \leq &\displaystyle \frac{2Me
  ^k|x|^{k+2}(\alpha_k|x|+\beta)^k}{(k+1)\big(r_0d_0(1-s)\big)^k}\int_0^1\frac{u^{k+1}}{\big((1-u^4)r_0+(s-u^4s^{\prime})d_0\big)^2} \, du.
\end{array} $$ 
Or, nous avons la majoration suivante~:
$$ \displaystyle
\int_0^1\frac{u^{k+1}}{\big((1-u^4)r_0+(s-u^4s^{\prime}d_0\big)^2} du \leq
\int_0^1 \frac{u}{\big((1-u^2)r_0+(s-u^2s^{\prime})d_0\big)^2} du \leq \frac{1}{2r_0d_0(s-s^{\prime})}.$$
Par suite, nous en d\'eduisons que~:
\begin{equation}\label{IneqI2}
I_2 \leq
\frac{Me^k(\alpha_k|x|+\beta)^k|x|^{k+2}}{(k+1)\big((r_0d_0(1-s)\big)^k}\Big(\frac{1}{r_0d_0(s-s^{\prime})}\Big).
\end{equation}
 
Par ailleurs, nous avons \'egalement~:
$$ \begin{array}{lll}
I_3 & = &\displaystyle 2\int_0^1 du \int_0^{u|x|}\frac{\|\delta_{k} (t)
  \|_{U_{s}}}{(1-u^4)r_0+(s-u^4s^{\prime})d_0} \, dt \\
\\   
& \leq & \displaystyle
\frac{2Me
  ^k(\alpha_k|x|+\beta)^k}{\big(r_0d_0(1-s)\big)^k}\int_0^1 du
\int_0^{u|x|} \frac{t^{k}}{(1-u^4)r_0+(s-u^4s^{\prime})d_0} \, dt\\
\\
&  \leq &\displaystyle \frac{2Me
  ^k(\alpha_k|x|+\beta)^k|x|^{k+1}}{(k+1)\big(r_0d_0(1-s)\big)^k}\int_0^1\frac{u^{k+1}}{(1-u^4)r_0+(s-u^4s^{\prime})d_0} \, du.
\end{array} $$ 
Or, nous avons les majorations suivantes~:
$$ \displaystyle
\int_0^1\frac{u^{k+1}}{(1-u^4)r_0+(s-u^4s^{\prime})d_0} du \leq
\int_0^1 \frac{u}{(1-u^2)r_0+(s-u^2s^{\prime})d_0} du,$$
d'o\`u
$$ \displaystyle \int_0^1\frac{u^{k+1}}{(1-u^4)r_0+(s-u^4s^{\prime})d_0} du \leq
\frac{1}{2(r_0+s^{\prime}d_0)}\ln\big(\frac{r_0+sd_0}{d_0(s-s^{\prime})}\big)
\leq \frac{r_1}{2r_0d_0(s-s^{\prime})}.$$
Par suite, nous en d\'eduisons que~:
\begin{equation}\label{IneqI3}
I_3 \leq
\frac{Me^k(\alpha_k|x|+\beta)^k|x|^{k+1}}{(k+1)\big(r_0d_0(1-s)\big)^k}\Big(\frac{r_1}{r_0d_0(s-s^{\prime})}\Big).
\end{equation}

Enfin, nous avons~:
$$ \begin{array}{lll}
I_4 & = &\displaystyle \int_0^1 du \int_0^{u|x|}t\|F\|_{D(0,r_1)}\|\delta_{k} (t)\|_{U_{s}} \, dt \\
\\   
& \leq & \displaystyle
\frac{Me
  ^k(\alpha_k|x|+\beta)^k \|F\|_{D(0,r_1)}}{\big(r_0d_0(1-s)\big)^k}\int_0^1 du
\int_0^{u|x|} t^{k+1} \, dt\\
\\
&  \leq &\displaystyle \frac{Me
  ^k|x|^{k+2}(\alpha_k|x|+\beta)^k \|F\|_{D(0,r_1)}}{(k+1)\big(r_0d_0(1-s)\big)^k}
\end{array} $$ 
d'o\`u~:
\begin{equation}\label{IneqI4}
I_4 \leq
\frac{Me^k(\alpha_k|x|+\beta)^k|x|^{k+2}}{(k+1)\big(r_0d_0(1-s)\big)^k}\Big(\frac{r_0d_0(s-s^{\prime})\|F\|_{D(0,r_1)}}{r_0d_0(s-s^{\prime})}\Big).
\end{equation}
 
Au final, en utilisant les majorations (\ref{IneqI1}), (\ref{IneqI2}),
(\ref{IneqI3}) et (\ref{IneqI4}), nous obtenons~:
\begin{equation}\label{IneqFin}
\|\delta_{k+1}\|_{U_{s^{\prime}}} \leq
\frac{Me^k(\alpha_k|x|+\beta)^k|x|^{k+1}}{\big(r_0d_0(1-s)\big)^k(k+1)}\Big[\big(1+r_0d_0(s-s^{\prime})\|F\|_{D(0,r_1)}\big)|x|+\beta\Big]\frac{1}{r_0d_0(s-s^{\prime})}.
\end{equation}

En choisissant dans (\ref{IneqFin})~:
$ \displaystyle s=s^{\prime}+\frac{1-s^{\prime}}{k+1}$, i.e
$\displaystyle 1-s=\frac{k}{k+1}(1-s^{\prime})$, nous obtenons~:
\begin{equation}
 \|\delta_{k+1}\|_{U_{s^{\prime}}} \leq
\frac{Me^k(\alpha_{k+1}|x|+\beta)^k|x|^{k+1}}{\big(r_0d_0(1-s^{\prime})\big)^k}\Big(1+\frac{1}{k}\Big)^k\frac{1}{r_0d_0(1-s^{\prime})}.
\end{equation}

Mais, $\forall k \in \mathbb{N}$, $\displaystyle
\big(1+\frac{1}{k}\big)^k \leq e$, donc nous en d\'eduisons finalement
que~:
\begin{equation}
 \displaystyle \|\delta_{k+1}\|_{U_{s^{\prime}}} \leq
M\Big[\frac{e|x|(\alpha_{k+1}|x|+\beta)}{r_0d_0(1-s^{\prime})}\Big]^{k+1},
\end{equation} 
ce qui ach\`eve la r\'ecurrence.
\end{proof}

\item Remarquons alors que 
$$ (\alpha_k|x|+\beta)^k=(|x|+3r_1)^k\bigg(1+\frac{r_0d_0(1-s)\|F\|_{D(0,r_1)}|x|}{(|x|+3r_1)k}\bigg)^k.$$
Par suite, par le lemme \ref{lemfon}, nous en d\'eduisons que la
s\'erie majorante de $\displaystyle \sum_{k \geq 0} \delta_k(x)$ converge d\`es que $$\displaystyle \frac{e|x|(|x|+3r_1)}{r_0d_0(1-s)}<1,$$
c'est-\`a-dire pour $$\displaystyle 0\leq|x|\leq \frac{3r_1}{2e}\bigg(-1+\sqrt{1+\frac{4r_0d_0(1-s)}{9er_1^2}}\bigg).$$
Par cons\'equent, la s\'erie $\displaystyle \sum_{k \geq 0}
\delta_k(x)$ converge absolument dans  $H(\overline{U_s},\mathbb{C})$
(pour tout $0 \leq s <1$) et
  uniform\'ement en $x \in K$, o\`u $K$ est un compact quelconque du
  disque ouvert $\displaystyle |x| < \min \bigg\{\frac{3r_1}{2e}\bigg(-1+\sqrt{1+\frac{4r_0d_0(1-s)}{9er_1^2}}\bigg),
  R\bigg\}$. Par construction, sa somme $\theta (x)$ satisfait l'\'equation~:
\begin{equation}\label{pre107}
\theta (x) =  \theta_0(x) + \int_0^1 du
\int_0^{ux} dt \,  L{(u,t)} \theta (t),
\end{equation}
de sorte que la fonction holomorphe $\varphi(z,x) := \theta(x)(z)$ est
solution de l'\'equation int\'egrale (\ref{le0ter}). En sp\'ecialisant
le r\'esultat pour $s=0$, nous obtenons le th\'eor\`eme par le lemme \ref{intrepre}. 
\end{enumerate}
\end{proof}

Nous d\'eduisons facilement du th\'eor\`eme \ref{prop22} le r\'esultat suivant~:
\begin{cor}\label{cor37}
Dans le th\'eor\`eme \ref{prop22}, si $F$ et $h$ sont des fonctions
enti\`eres de $z$, alors
$\psi(z,x)$ s'\'etend analytiquement \`a $\mathbb{C}^2$.
\end{cor}

\subsubsection{Construction explicite}

Nous revenons maintenant \`a la fonction  $\Psi(z, \z)$ associ\'ee \`a $\displaystyle
\psi(z,x)$ par (\ref{eq18}). Du lemme \ref{lem2}, 
du th\'eor\`eme \ref{prop22} et de son corollaire \ref{cor37}, nous
d\'eduisons le r\'esultat suivant~:

\begin{prop}\label{bilant}
Supposons que $F$ et $h$ soient des fonctions holomorphes au voisinage
de l'origine. Alors il existe une unique fonction holomorphe
$\Psi(z,\z)$ au voisinage de $(0,0)$
satisfaisant les conditions~:
\begin{equation}
 \left\{
\begin{array}{l}
\displaystyle \Psi(z,\z) \big|_{\displaystyle z=\z^2} =1  \\
\\
\displaystyle \bigg(-\frac{1}{2\z}\,\frac{\partial \Psi}{\partial \z}(z,\z)\bigg)
\big|_{\displaystyle z=\z^2} = h(z).
\end{array}
 \right. 
\end{equation}
et telle que  $\displaystyle \tpsi(z,\z) := \frac{\Psi(z, \z)}{z-\z^2}$ soit solution de
l'EDP lin\'eaire (\ref{eq9}).\\
De plus, si $F$ et $h$ sont des fonctions enti\`eres, alors $\Psi(z,\z)$ s'\'etend
analytiquement \`a $\mathbb{C}^2$.\\
\end{prop}

Par suite, nous avons facilement~:

\begin{prop}\label{prop251}
Supposons que $F$ et $h$ soient des fonctions holomophes au voisinage
de l'origine. Alors la fonction $\stackrel{\vee}{\Phi}(z,\xi)$ d\'efinie par~:
\begin{equation}\label{eq27}
 \stackrel{\vee}{\Phi}(z,\xi)_{\displaystyle | \xi = S(z,\z)} :=
 \frac{\Psi(z,\z)}{z-\z^2}
\end{equation}
avec $\Psi$ comme dans la proposition \ref{bilant}, est solution de
(\ref{eq8}) et est un majeur d'une microfonction confluente en $(0,0)$
  \`a support singulier la courbe alg\'ebrique $\mathcal{C} = \{(z,
  \xi),\,\, 9\xi^2 = 4z^3 \}$ (cf. d\'efinition
\ref{15}). \\
Lorsque $F$ et $h$ sont des fonctions enti\`eres, alors
$\stackrel{\vee}{\Phi}(z,\xi)$ est 
un majeur d'une microfonction confluente {\em r\'esurgente} en $(0,0)$
  \`a support singulier la courbe alg\'ebrique $\mathcal{C}$.
\end{prop}

L'existence d'un tel majeur va nous permettre de construire in fine
les fonctions confluentes recherch\'ees en utilisant la
repr\'esentation int\'egrale (\ref{eq7}).

\begin{figure}[thp]
\begin{center}
\begin{tabular}{|ccc|}
\hline
 \includegraphics[width=2.1in]{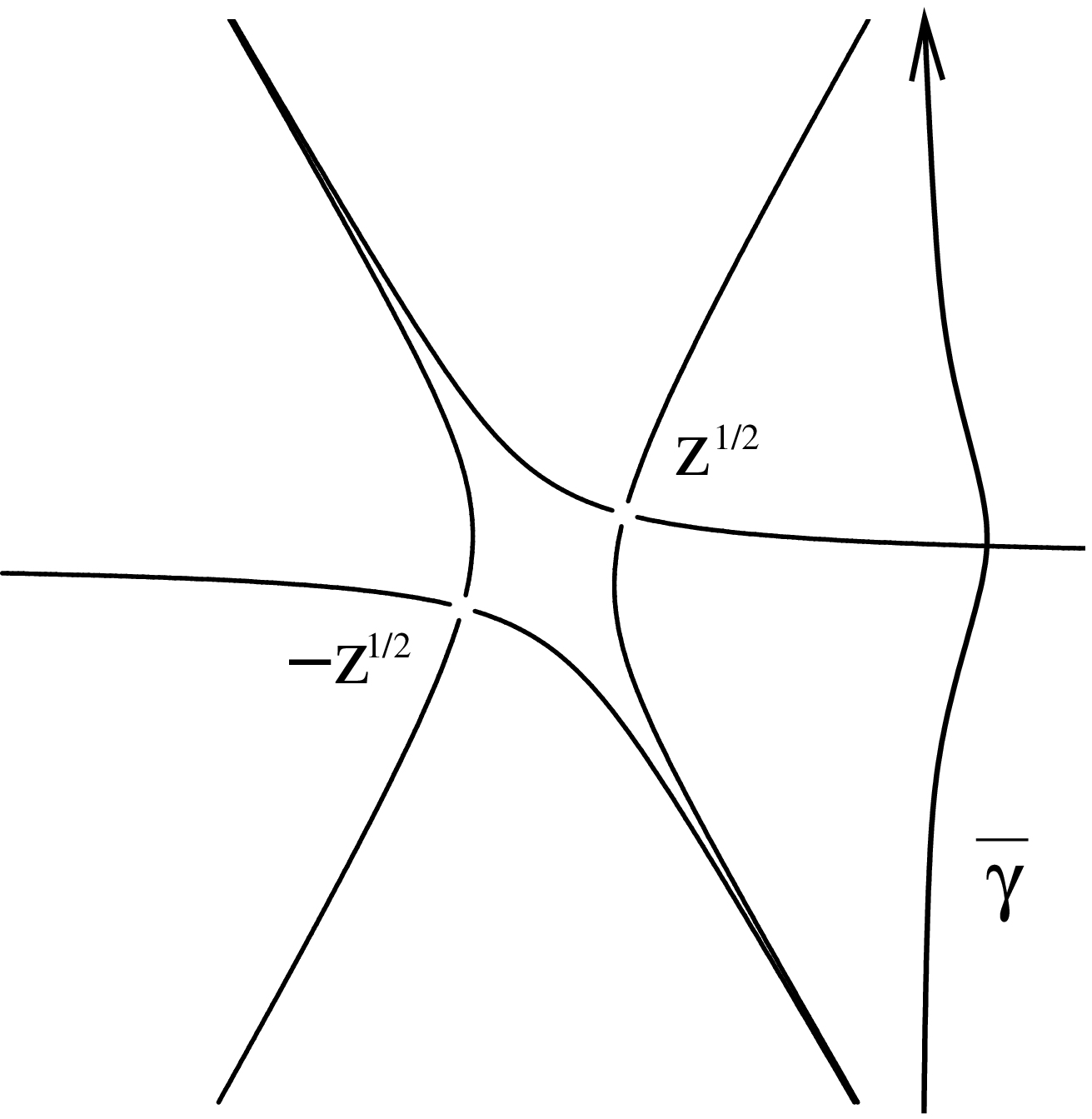} & \hspace{10mm} &
\includegraphics[width=1.7in]{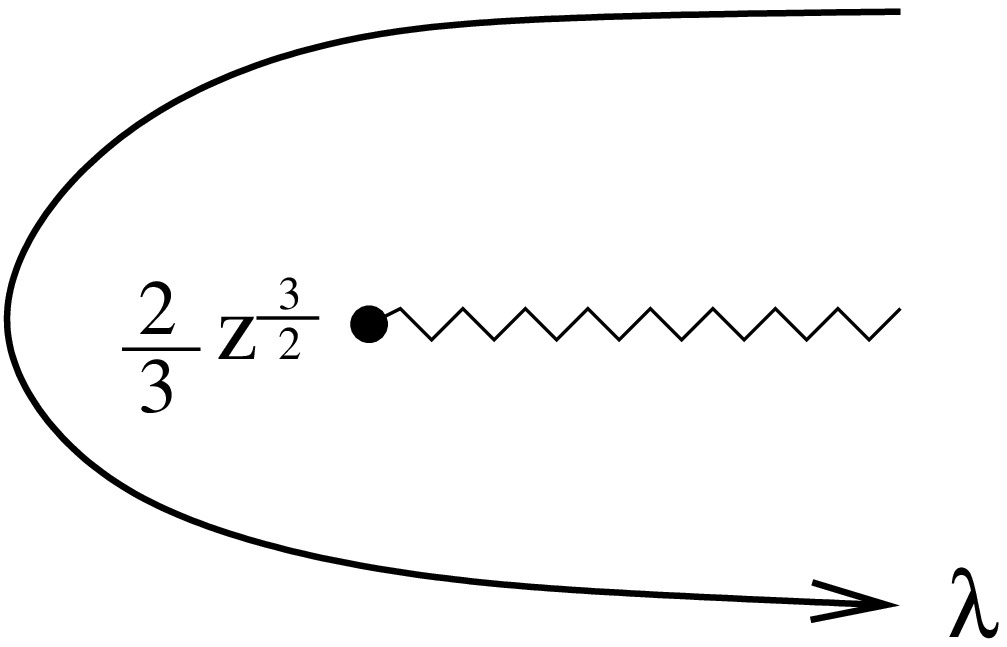} \\
Fig. \ref{fig:Stokes66}.1a  &  &  Fig. \ref{fig:Stokes66}.1b\\ 
\hline
\includegraphics[width=2.1in]{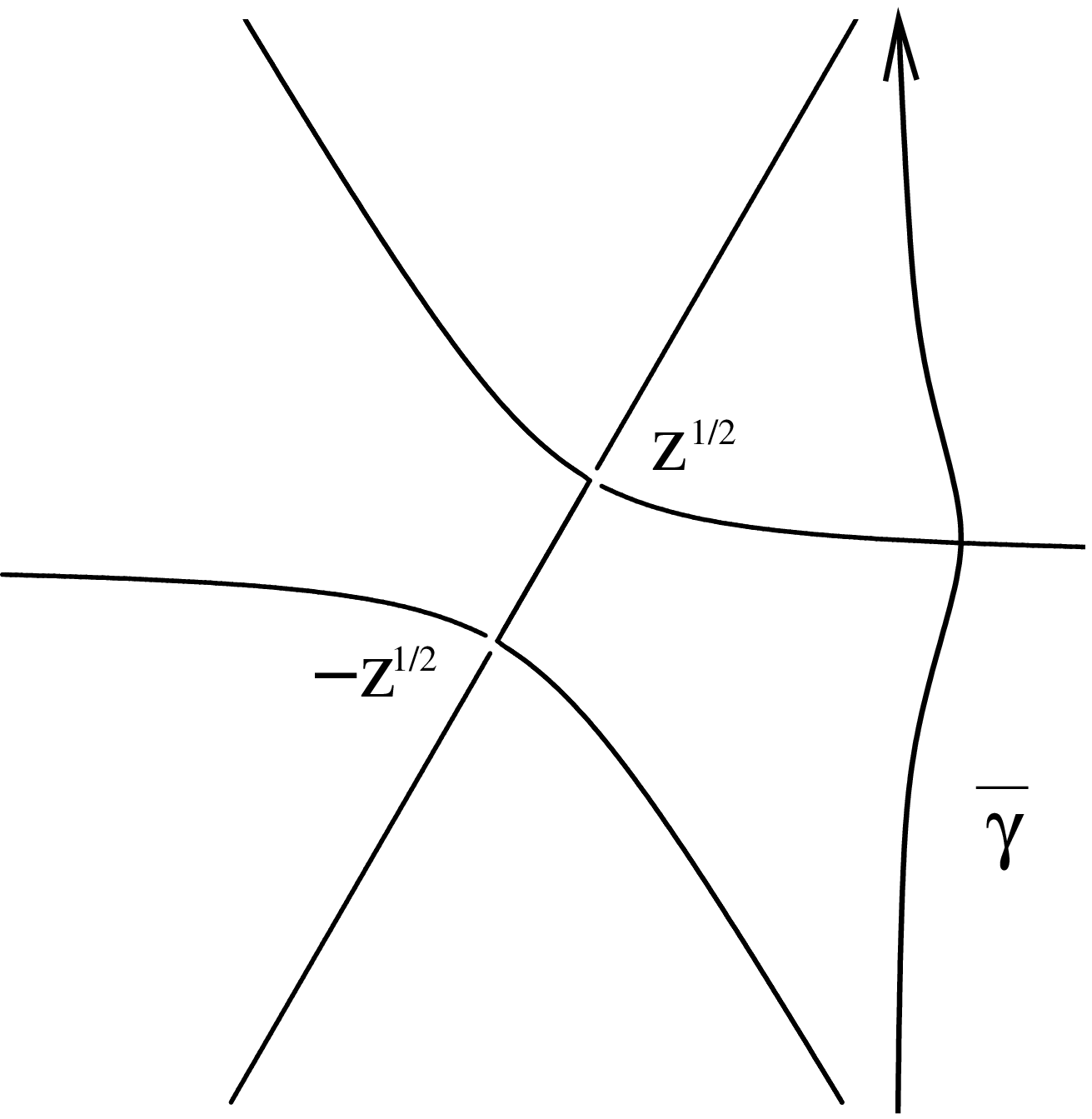} & &
\includegraphics[width=1.7in]{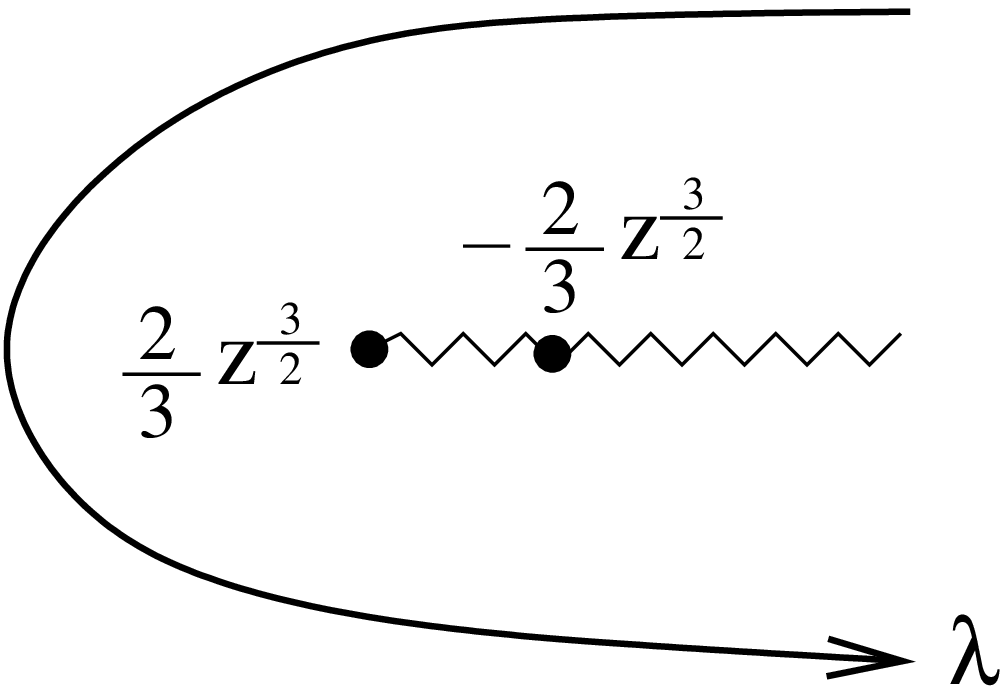}\\
Fig. \ref{fig:Stokes66}.2a  &  &  Fig. \ref{fig:Stokes66}.2b\\ 
\hline
\includegraphics[width=2.1in]{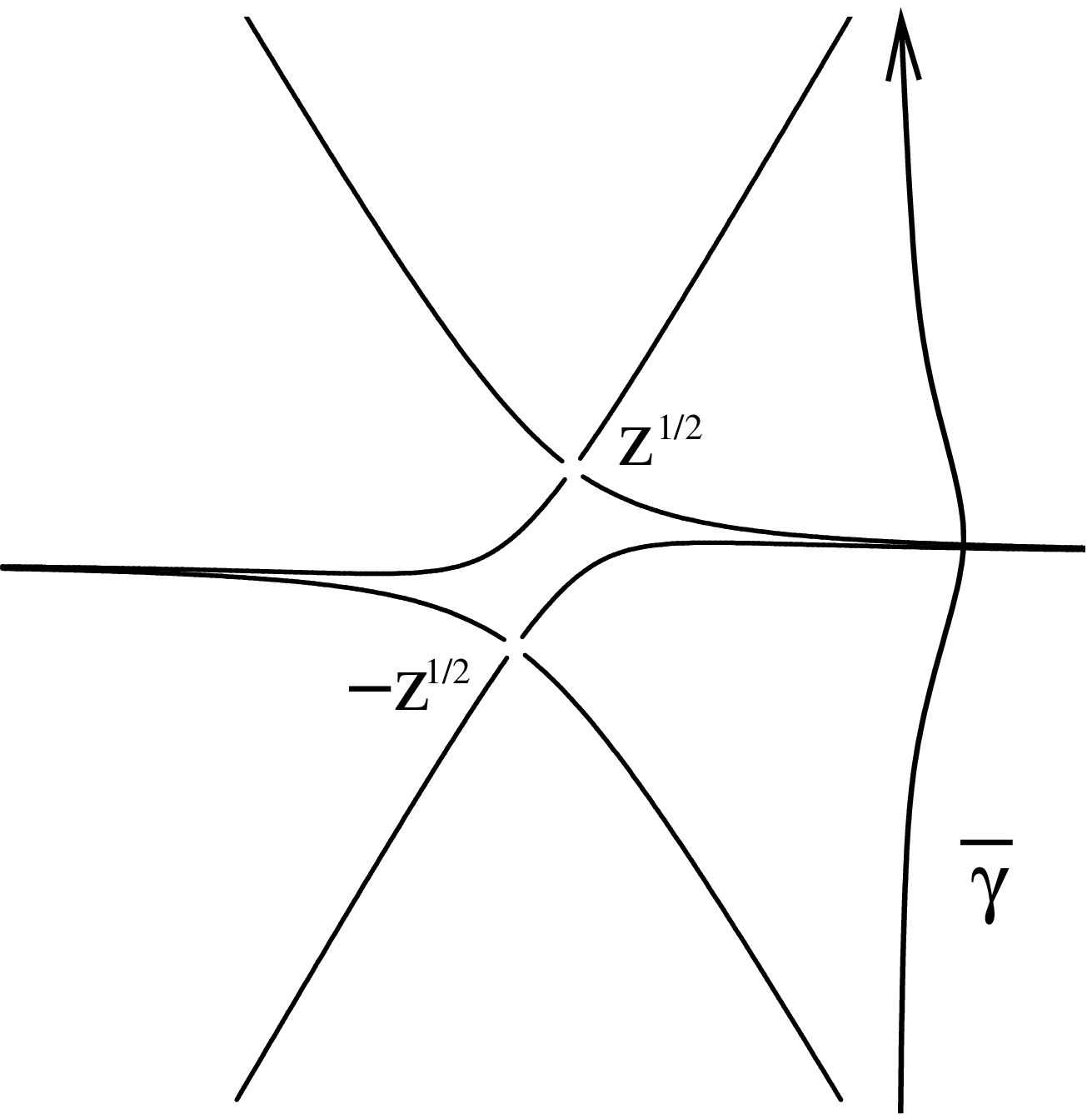} & &
\includegraphics[width=1.7in]{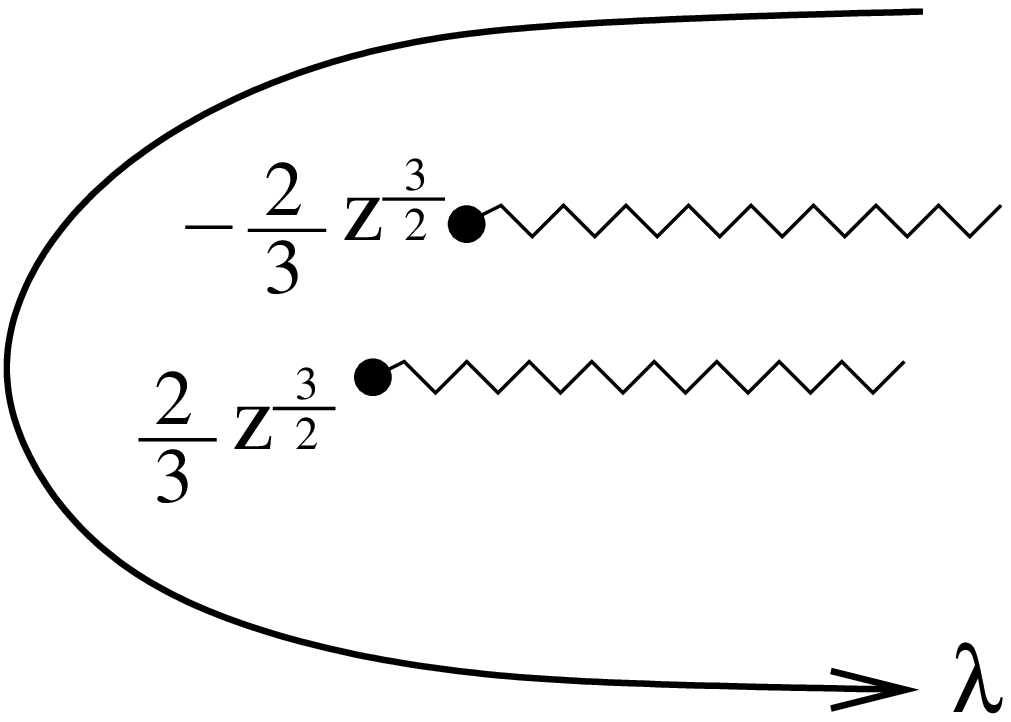}\\
Fig. \ref{fig:Stokes66}.3a  &  &  Fig. \ref{fig:Stokes66}.3yb\\ 
\hline
\end{tabular}
\caption{Sur les figures de gauche les chemins de plus grande pente
et le chemin $\overline \gamma$, 
  et sur les figures de droite son image $\overline \lambda$ par la transformation $\z \mapsto \xi =
  S(z,\z)$ (les lignes ondul\'ees sont coup\'ees). Fig. \ref{fig:Stokes66}.1 pour $z \in S_1$, Fig. \ref{fig:Stokes66}.2
  pour $z$ sur $L_1$, Fig. \ref{fig:Stokes66}.3 pour $z \in S_2$ (voir
  figure \ref{fig:Stokes5}). 
\label{fig:Stokes66}}
\end{center}
\end{figure}

\begin{prop}\label{fonctconfl}
Consid\'erons la repr\'esentation int\'egrale
\begin{equation}\label{eq7bis}
\displaystyle \Phi(z,\varepsilon) =  
\int_\gamma  e^{-\frac{1}{\varepsilon} S(z, \widehat{z})} \Psi(z,\widehat{z})  \, d\widehat{z}
\end{equation}
avec $\Psi$ comme dans la proposition \ref{bilant} (ou d'une mani\`ere
\'equivalente la repr\'esentation int\'egrale
\begin{equation}\label{eq26}
\displaystyle \Phi(z,\varepsilon) =  
\int_{\lambda}  e^{-\frac{\xi}{\varepsilon}} \stackrel{\vee}{\Phi}(z,\xi) \, d\xi.
\end{equation}
avec $\stackrel{\vee}{\Phi}(z,\xi)$ comme dans la proposition \ref{prop251}).\\
Notons $\overline{\gamma}$ le chemin $\gamma$ tronqu\'e comme dans la figure
\ref{fig:Stokes66} (et $\overline \lambda$ son image par la transformation $\z \mapsto \xi = 
  S(z,\z)$ pour la repr\'esentation (\ref{eq26})).\\
Alors, si $F$ et $h$ sont holomorphes au voisinage de l'origine
(respectivement enti\`eres), alors les repr\'esentations int\'egrales
\begin{equation}\label{eq260}
\displaystyle  \Phi(z,\varepsilon) =  
\int_{\overline \gamma}  e^{-\frac{1}{\varepsilon} S(z, \z)} \Psi(z,\z)  \, d\z,
\end{equation}
et \begin{equation}\label{eq26bis}
\displaystyle \Phi(z,\varepsilon) =  
\int_{\overline  \lambda}  e^{-\frac{\xi}{\varepsilon}} \stackrel{\vee}{\Phi}(z,\xi) \, d\xi
\end{equation}
repr\'esentent une fonction confluente $\bold \Phi (z, \varepsilon)$ (respectivement une fonction
confluente r\'esurgente) \`a support singulier la courbe alg\'ebrique
$\mathcal{C}$ (au sens de la d\'efinition \ref{16}). 
\end{prop}

\begin{proof}
\begin{enumerate}
\item Dans le cas o\`u $F$ et $h$ sont holomorphes au voisinage de
  l'origine,  en suivant la proposition \ref{bilant}, nous savons que  $\Psi(z, \z)$ est holomorphe
dans un voisinage de l'origine dans $\mathbb{C}^2$, disons pour $\displaystyle (z,\z) \in
D(0,\frac{r^2}{4})\times  D(0,r)$ avec $r>0$ assez petit, o\`u
$D(0,r)$ d\'esigne le disque ouvert de rayon $r$ centr\'e en $0$. Par cons\'equent, l'int\'egrale
(\ref{eq7}) est bien d\'efinie pourvu que nous tronquions le chemin
d'int\'egration $\gamma$ qui est alors not\'e par
$\overline \gamma$, comme sur la figure \ref{fig:Stokes66}.\\
Alors, dans la repr\'esentation int\'egrale 
\begin{equation}\label{eq260bis}
\displaystyle  \Phi(z,\varepsilon) =  
\int_{\overline \gamma}  e^{-\frac{1}{\varepsilon} S(z, \z)} \Psi(z,\z)  \, d\z,
\end{equation}
en faisant le changement de variable $\xi = S(z, \z)$, 
nous avons l'int\'egrale correspondante, o\`u le chemin $\overline  \lambda$ est dessin\'e
sur la figure \ref{fig:Stokes66}.
\begin{equation}\label{eq26bisbis}
\displaystyle \Phi(z,\varepsilon) =  
\int_{\overline  \lambda}  e^{-\frac{\xi}{\varepsilon}} \stackrel{\vee}{\Phi}(z,\xi) \, d\xi.
\end{equation}
Une cons\'equence de la proposition \ref{prop251} est que
l'int\'egrale de  Laplace  (\ref{eq26bisbis})  repr\'esente une fonction confluente $\bold \Phi (z, \varepsilon)$   \`a support singulier dans $\mathcal{C}$, au sens de la d\'efinition \ref{16}.\\

\item Lorsque $F$ et $h$ sont des fonctions  enti\`eres, puisque par la proposition
\ref{bilant} ({\em resp} proposition \ref{prop251}) $\Psi(z,\z)$
({\em resp} $\stackrel{\vee}{\Phi}(z,\xi)$) 
s'\'etend analytiquement dans tout $\mathbb{C}^2$ ({\em resp} s'\'etend
comme un majeur d'une microfonction confluente r\'esurgente), l'int\'egrale tronqu\'ee
 (\ref{eq260bis}) ({\em resp.} (\ref{eq26bisbis})) a encore un sens
pour toute tronquature, et nous pouvons interpr\'eter les
repr\'esentations int\'egrales (\ref{eq26bisbis}) et (\ref{eq260bis})
comme une {\em pr\'esomme de Borel} \cite{DP99, CNP2}, d\'efinissant ainsi
  fonction confluente r\'esurgente $\bold \Phi (z, \varepsilon)$ \`a support
  singulier dans 
$\mathcal{C}$ (cf. remarque \ref{16bis}).
\end{enumerate}
\end{proof}

\subsection{D\'ecomposition et cons\'equences}

La d\'ecomposition locale  ({\em resp.} d\'ecomposition) de la fonction
confluente $\bold \Phi (z, \varepsilon)$ ({\em resp} fonction confluente
r\'esurgente $\bold \Phi (z, \varepsilon)$) de la proposition \ref{fonctconfl}
peut se d\'eduire de la repr\'esentation int\'egrale 
 (\ref{eq260}) par la m\'ethode du col.
Nous d\'ecrivons ce que nous obtenons pour un germe de secteurs de Stokes ({\em
  resp.} secteurs de Stokes) dans la figure \ref{fig:Stokes5}.a.

\subsubsection{D\'ecomposition dans $S_1$}

\begin{lem}
Pour $z$ dans le germe de secteurs de Stokes ({\em
  resp.} secteur de Stokes) $S_1$, la d\'ecomposition locale
({\em resp.} d\'ecomposition) de la fonction confluente 
({\em  resp.} fonction confluente r\'esurgente) $\bold \Phi (z, \varepsilon)$
  induit un d\'eveloppement BKW formel unique~:
\begin{equation}\label{eqa28}
\begin{array}{ccc}
\bold \Phi (z, \varepsilon) & \stackrel{\displaystyle \sigma_{S_1}}{\displaystyle
  \longrightarrow} & i\sqrt{\pi\varepsilon}\Phi_{bkw}^+(z,\varepsilon),
\end{array}
\end{equation}
\end{lem} 

\begin{proof}
Pour la repr\'esentation int\'egrale (\ref{eq260}), $z$ \'etant dans
le (germe de) secteur de Stokes $S_1$, cela correspond \`a la situation
d\'ecrite sur la Fig. \ref{fig:Stokes66}.1a. En d\'eformant le chemin
d'integration $\overline \gamma$ sous le flot $\displaystyle \nabla \left( \Re \Big(
  \frac{S(z,\z)}{\varepsilon} \Big)\right)$ (les extr\'emit\'es
$\overline \gamma$ restant fix\'ees), nous voyons que seul le point col $\z =
\sqrt{z}$ a une contribution non triviale \`a la d\'ecomposition. Ceci
donne la formule (\ref{eqa28}). Notons que le d\'eveloppement BKW
formel ainsi obtenu est une solution formelle de l'\'equation
(\ref{eq4}) puisque le majeur $\stackrel{\vee}{\Phi}(z,\xi)$ est une solution de (\ref{eq8})
(cf. Prop. \ref{prop251}).
\end{proof}

De la m\^eme mani\`ere, nous pouvons montrer que, pour $z$ sur la ligne
de Stokes $L_1$, un ph\'enom\`ene de Stokes se produit
  (voir Fig. \ref{fig:Stokes66}.2) de sorte que, pour $z$ dans le germe
  de secteur de Stokes ({\em
  resp.} secteur de Stokes) $S_2$  (voir Fig. \ref{fig:Stokes66}.3a), la d\'ecomposition locale
 ({\em resp.} d\'ecomposition) de la fonction confluente  $\bold \Phi
 (z, \varepsilon)$ induit maintenant une somme de deux d\'eveloppements BKW formels~:
\begin{equation}\label{eqa29}
\begin{array}{ccc}
\bold \Phi (z, \varepsilon) & \stackrel{\displaystyle \sigma_{S_2}}{\displaystyle
  \longrightarrow} & i\sqrt{\pi\varepsilon}\Phi_{bkw}^+(z,\varepsilon) - \sqrt{\pi\varepsilon}
  \Phi_{bkw}^-(z,\varepsilon),
\end{array}
\end{equation}
o\`u $\Phi_{bkw}^-(z, \varepsilon) = \Phi_{bkw}^+(z,-\varepsilon)$.

R\'eciproquement, consid\'erons la solution BKW \'el\'ementaire 
$\Phi_{bkw}(z,\varepsilon)$ de l'\'equation (\ref{eq4}). Par  (\ref{eqa28}),
\`a un facteur
$\displaystyle i\sqrt{\pi\varepsilon}$ pr\`es, une d\'etermination de 
cette solution BKW \'el\'ementaire appara\^it comme la d\'ecomposition
locale dans un germe de secteur de Stokes d'une fonction confluente.
Plus g\'en\'eralement, toute d\'etermination de 
$\Phi_{bkw}(z,\varepsilon)$ dans n'importe quel germe de secteur de  Stokes
peut \^etre vue comme la d\'ecomposition locale dans ce germe de
secteur de Stokes d'une fonction confluente
(\`a un facteur $\displaystyle c\sqrt{\pi\varepsilon}$, $c \in \{\pm 1, \pm i\}$ pr\`es): dans la
repr\'esentation int\'egrale (\ref{eq260}), cela en d\'ecoule
simplement en choisissant un chemin d'int\'egration tronqu\'e
convenable $\overline \gamma$. En outre, puisque
toute notre analyse peut \^etre reconduite en choisissant une autre direction
$\alpha$ que $0$, nous obtenons le r\'esultat suivant~:

\begin{thm}\label{propdecomp}
Lorsque $F$ est holomorphe au voisinage de l'origine (respectivement enti\`ere), il existe une
famille de solutions BKW \'el\'ementaires (respectivement
\'el\'ementaires r\'esurgentes) de {\em type Airy local}
(respectivement de {\em type Airy})
$\Phi_{bkw}(z,\varepsilon)$ de l'\'equation (\ref{eq4}) (au
sens de la d\'efinition \ref{deftypeAiry}).
\end{thm}

\subsubsection{Lien avec le mod\`ele d'Airy}

Nous d\'eduisons alors du th\'eor\`eme \ref{propdecomp} pr\'ec\'edent,
en appliquant simplement un th\'eor\`eme de Jidoumou \cite{Jidoumou}~:

\begin{thm}\label{decompose}
Si $F$ est une fonction holomorphe au voisinage de l'origine
(respectivement enti\`ere), et si nous notons
$\Phi_{bkw}(z,\varepsilon)$ une solution BKW \'el\'ementaire
(respectivement \'el\'ementaire r\'esurgente)
donn\'ee par le th\'eor\`eme \ref{propdecomp}, alors pour $z\neq 0$
dans un voisinage de $0$ (respectivement $z \in \mathbb{C} \backslash \{0\}$), nous avons la d\'ecomposition unique suivante~:
$$\Phi_{bkw}(z,\varepsilon) = a(z, \varepsilon) A_{bkw}(z,\varepsilon) + b(z,\varepsilon)
\,\varepsilon\frac{\partial A_{bkw}}{\partial z}(z, \varepsilon),$$
o\`u $A_{bkw}(z,\varepsilon)$ d\'esigne le symbole BKW d'Airy
(\ref{wkb10}), tandis que
$a$ et $b$ (qui d\'ependent de $\Phi_{bkw}$) sont des constantes
locales de r\'esurgence (respectivement constantes de r\'esurgence), $a$ inversible et $b$ petite.
\end{thm}

Notons qu'ici, ``une constante (locale) de r\'esurgence'' signifie la
chose suivante~:

\begin{defn}
Une {\em constante locale de r\'esurgence} (respectivement {\em
  constante de r\'esurgence}) $c(Z, \eta)$ est un
d\'eveloppement BKW formel $\displaystyle c(Z,\varepsilon) = \sum_{n \geq 0} c_n(Z)
\varepsilon ^{n}$ tel que son mineur $\displaystyle \sum_{n \geq 1} c_n(Z)
\frac{\xi^{n-1}}{\Gamma(n)}$ d\'efinit un germe de fonctions
holomorphes en
$(Z,\xi) =(0,0) \in \mathbb{C}^m \times \mathbb{C}$, $m \geq 1$
(respectivement une fonction holomorphe sur $\mathbb{C}^m \times \mathbb{C}$, $m \geq 1$).
\end{defn} 

\section{Applications}\label{sec5}

\subsection{Un th\'eor\`eme local de r\'eduction}
En reproduisant le raisonnement de Pham (\cite{Ph00} \S 2.4), nous
d\'eduisons du th\'eor\`eme \ref{decompose} le th\'eor\`eme suivant~:

\begin{thm}\label{theo42}
Supposons que $F(z)$ dans (\ref{eq4}) est holomorphe au voisinage de l'origine
({\em resp.} une fonction enti\`ere). Alors
il existe une constante locale de r\'esurgence ({\em resp.} une
constante de r\'esurgence)  $s(z,\varepsilon)$ telle que, sous l'action
de la transformation
\begin{equation}\label{eq42bis}
 \left\{
 \begin{array}{l}
\displaystyle s(z, \varepsilon) = \sum_{k\geq 0} s_k(z)\varepsilon^{k}, \hspace{5mm} s_0(z) = z \\
\\
\displaystyle \Phi(z, \varepsilon)=\Big( \frac{\partial s}{\partial z} \Big)^{-\frac{1}{2}}
y \big( s(z, \varepsilon) , \varepsilon\big),
\end{array}
 \right.
\end{equation}
l'\'equation (\ref{eq4}) devient l'\'equation (\ref{eq3}),
pour $z$ ({\em resp.} $s$) au voisinage de l'origine. De plus, sous
l'action (\ref{eq42bis}), une solution BKW
\'el\'ementaire de (\ref{eq3}) est
transform\'ee en une solution BKW \'el\'ementaire de (\ref{eq4}).
\end{thm}

\subsection{Applications pour l'\'equation de Schr\"odinger}

Nous nous focalisons maintenant sur l'\'equation de Schr\"odinger~:
\begin{equation}\label{eq1}
 \varepsilon^2 \frac{d^2 Y}{d q^2}=V(q)Y
\end{equation}
avec 
\begin{equation}\label{eq2}
V(q) =  q+\sum_{n=2}^{+\infty} v_n q^n 
\end{equation}
analytique au voisinage de l'origine admettant $0$ comme z\'ero
simple.\\
Cette hypoth\`ese sur $V$ signifie que $q=0$ est un point tournant
simple pour les solutions formelles BKW.
Nous rappelons que ces solutions BKW sont des combinaisons lin\'eaires
de solutions formelles BKW \'el\'ementaires de (\ref{eq1}) de la
forme :
\begin{equation}\label{eq2bis}
Y_{bkw} (q, \varepsilon) = e^{\displaystyle  -\frac{1}{\varepsilon} \int^q \sqrt{ V(t) } dt}
\left( \sum_{k=0}^{+\infty} Y_k(q) \varepsilon^{k}  \right).
\end{equation}
Nous rappelons \'egalement que ces solutions formelles BKW
\'el\'ementaires (d\'efinies localement en $q$) sont d\'efinies de
mani\`ere unique \`a normalisation pr\`es, i.e \`a multiplication
pr\`es par un d\'eveloppement formel inversible de la forme
$\displaystyle  e^{\frac{c}{\varepsilon}}\sum_{k \geq 0}  c_k \varepsilon^{k}$, $c, c_k \in
\mathbb{C}$, $c_0 \neq 0$.\\ 
Nous voulons traduire l'analyse BKW que nous avons faite pour
l'\'equation (\ref{eq4}) en une analyse analogue pour l'\'equation
(\ref{eq1}).\\
En suivant \cite{Aoki91} et \cite{Ph00}, le premier pas naturel afin
d'obtenir notre th\'eor\`eme de r\'eduction est de ``redresser'' la
g\'eom\'etrie au voisinage de l'origine via un changement de variable $q
\leftrightarrow z$ qui transforme la forme diff\'erentielle
$\displaystyle \sqrt{V(q)} dq$ en $\displaystyle \sqrt{z} dz$
(l'application cotangente associ\'ee transforme l'\'equation de la
sous-vari\'et\'e Lagrangienne $P^2-V(q)=0$ en $p^2-z=0$ ).\\
Ceci nous am\`ene \`a poser la transformation~:
\begin{equation}\label{Green}
 \left\{
 \begin{array}{l}
\displaystyle z(q)=\left( \frac{3}{2} \int_{0}^{q} V(t)^{\frac{1}{2}}
dt \right)^{\frac{2}{3}} \\
\\
\displaystyle Y(q, \varepsilon)=\Big( \frac{dz}{dq} \Big)^{-\frac{1}{2}}
\Phi (z, \varepsilon).
\end{array}
 \right.
\end{equation}
Ce changement de variable transforme (\ref{eq1}) en notre \'equation
``canonique''~:
$$ \displaystyle  \frac{d^2 \Phi}{d z^2} - \frac{z}{\varepsilon^2} \Phi = F(z) \Phi,$$ 
avec 
\begin{equation}\label{eq4bis}
\displaystyle F(z) =  \frac{z}{2V(q)} \{z,q\}_{\displaystyle |q=q(z)} ,
\end{equation}
o\`u  $\displaystyle \{z,q\}= \frac{z^{\prime \prime
    \prime}(q)}{z^{\prime}(q)}-\frac{3}{2}\left( \frac{z^{\prime
      \prime}(q)}{z^{\prime}(q)}\right)^2$ d\'esigne la d\'eriv\'ee Schwarzienne
de $z$ par rapport \`a $q$.\\ 

La fonction $F$ ainsi d\'efinie satisfait la propri\'et\'e suivante~:
\begin{lem}\label{regF}
Si $V(q)$ est holomorphe au voisinage de $0$, $V(q) \sim q$, 
alors $F(z)$ est holomorphe au voisinage de l'origine.
\end{lem}
 
\begin{proof}
Du d\'eveloppement en s\'erie de Taylor convergent (\ref{eq2}), nous
en d\'eduisons que $\displaystyle \frac{3}{2} \int_{0}^{q}  V(t)^{\frac{1}{2}}
dt -q^{\frac{3}{2}} \in q^{\frac{5}{2}}\mathbb{C}\{q\}$, de sorte que 
$\displaystyle z(q)-q \in q\mathbb{C}\{q\}$.\\
Le th\'eor\`eme de Lagrange nous permet d'obtenir la fonction inverse
$q(z)$ qui est elle aussi une fonction holomorphe au voisinage de $0$.
Etant donn\'e que  $\displaystyle 
\{z,q\}=-\{q,z\}\Big(\frac{dz}{dq}\Big)^2$, nous en d\'eduisons
facilement que $F(z)$ est holomorphe au voisinage de $0$ (et $\displaystyle F(z)
=\frac{3}{7}v_3-\frac{9}{35}v_2^2 + O(z)$, o\`u les $v_i$ sont
d\'efinis en (\ref{eq2})). 
\end{proof}

Lorsque $V(q)$ est une fonction enti\`ere (ou m\^eme une fonction
m\'eromorphe), un r\'esultat plus pr\'ecis peut \^etre obtenu en
utilisant les propri\'et\'es bien connues de la transformation $\displaystyle q \mapsto \int^{q} V(t)^{\frac{1}{2}}
dt$ en termes de transformation conforme (voir, par exemple, \cite{Jenkins, Sib75,
  Elsner}). La figure \ref{fig:Stokes1} illustre ce type de
r\'esultat~:
\begin{figure}[thp]
\begin{center}
\begin{tabular}{ccc}
\includegraphics[width=2.6in]{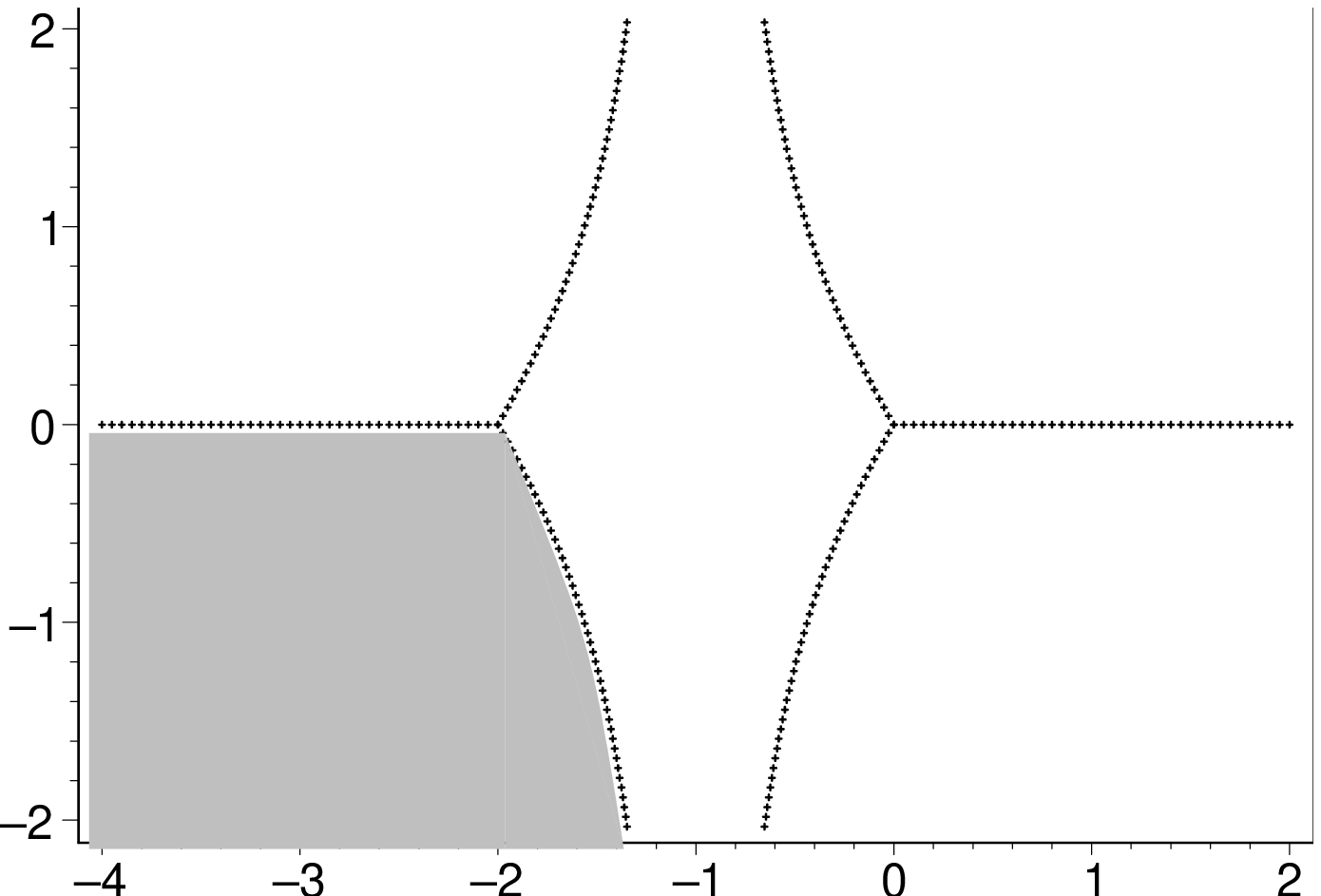}  & \hspace{5mm} &
\includegraphics[width=2in]{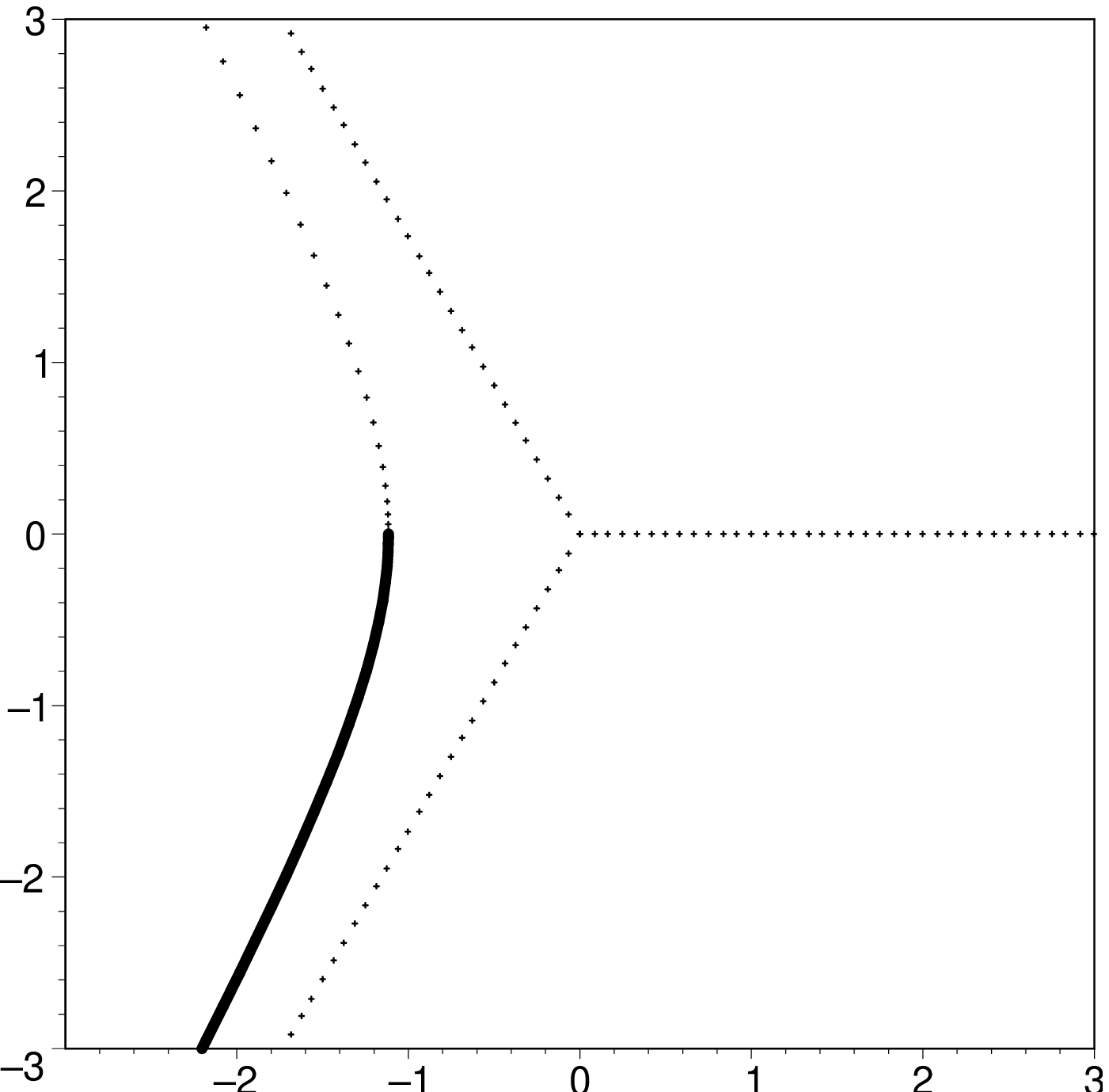} \\
Fig. \ref{fig:Stokes1}.a  &  &  Fig. \ref{fig:Stokes1}.b 
\end{tabular}
\caption{La transformation $q \rightarrow z$ (donn\'ee par (\ref{Green})) 
pour $\displaystyle V(q)
  = q + \frac{1}{2}q^2$. L'ouvert compl\'ementaire de la zone gris\'ee
 dans le $q$-plan (Fig. \ref{fig:Stokes1}.a) est envoy\'e
 conform\'ement sur le domaine coup\'e dans le $z$-plan
  (Fig. \ref{fig:Stokes1}.b. La coupure est la ligne pleine). Les
  lignes en pointill\'es sont les lignes de Stokes (relatives \`a la direction d'argument $0$). La
  fonction $F(z)$ est holomorphe dans le domaine coup\'e dessin\'e
  \`a la Fig. \ref{fig:Stokes1}.b.
\label{fig:Stokes1}}
\end{center}
\end{figure}

Par la transformation (\ref{Green}), nous d\'eduisons maintenant du
th\'eor\`eme \ref{theo42} le th\'eor\`eme suivant~:
\begin{thm}\label{theo1}
Il existe une constante locale de r\'esurgence $\displaystyle \big(s_k(q)\big)_{k\geq 0}$ telle que, sous l'action de la transformation
\begin{equation}\label{eq3bis}
 \left\{
 \begin{array}{l}
\displaystyle s(q, \varepsilon) = \sum_{k\geq 0} s_k(q)\varepsilon^{k} \\
\\
\displaystyle Y(q, \varepsilon  )=\Big( \frac{\partial s}{\partial q} \Big)^{-\frac{1}{2}}
y \big( s(q, \varepsilon) , \varepsilon\big),
\end{array}
 \right.
\end{equation}
l'\'equation (\ref{eq1}) est chang\'ee en l'\'equation (\ref{eq3}),
pour $q$ ({\em resp.} $s$) au voisinage de l'origine. De plus, sous
l'action de (\ref{eq3bis}), une solution BKW formelle (\ref{eq3}) est
transform\'ee en une solution BKW formelle de (\ref{eq1}).
\end{thm}

Autrement dit, nous avons montr\'e que, dans le cadre de l'analyse BKW,
l'\'equation (\ref{eq1}) se ram\`ene \`a l'\'equation (\ref{eq3}).

Au niveau formel, ce type de r\'esultat a d\'ej\`a \'et\'e \'etabli
dans un article de Silverstone \cite{Sil85}, et depuis d'autres r\'esultats
 plus pr\'ecis concernant les propri\'et\'es de la transformation $s$
 dans (\ref{eq3bis}) ont \'et\'e
 \'etablis~:
\begin{itemize}
 \item 
Dans \cite{Aoki91} (voir aussi \cite{KawaiTak94}), T. Aoki, T. Kawai et Y. Takei 
d\'emontrent le th\'eor\`eme \ref{theo1} : nous retrouvons ainsi leur
r\'esultat \`a la diff\'erence que nous n'avons \`a aucun moment
utilis\'e le calcul microdiff\'erentiel de Sato. Par la suite, ce r\'esultat a \'et\'e \'etendu dans \cite{Aoki93} au cas o\`u dans l'\'equation (\ref{eq1})
la fonction potentielle $V$ est une constante de r\'esurgence locale.
\item 
Dans \cite{Ph00}, F. Pham montre que
le d\'eveloppement $s(q, \varepsilon)$ est {\em r\'esurgent en $\varepsilon^{-1}$ 
\`a d\'ependance r\'eguli\`ere  en $q$}. Cependant, ce r\'esultat est bas\'e sur l'hypoth\`ese
qu'une base de solutions BKW r\'esurgentes de (\ref{eq1})
puisse \^etre d\'efinie, avec une d\'ependance r\'eguli\`ere en $q$
except\'e aux points tournants. 
\end{itemize}

\subsection{Extensions possibles}

Nous pouvons \'etendre facilement nos r\'esultats \`a l'\'equation~: 
\begin{equation}\label{eqq4}
\displaystyle  \frac{d^2 \Phi}{d z^2} - \frac{z}{\varepsilon^2} \Phi = F(z, \beta) \Phi,
\end{equation}
o\`u $F$ d\'epend holomorphiquement de $(z,\beta) \in \mathbb{C}^2$ au voisinage
de l'origine. Dans ce cas, les th\'eor\`emes \ref{propdecomp} et
\ref{decompose} deviennent~:

\begin{thm}\label{propdecompbis}
Il existe une famille de solutions BKW \'el\'ementaires
$\Phi_{bkw}(z, \beta,\varepsilon)$ de l'\'equation (\ref{eqq4}) qui sont
de type Airy local, \`a d\'ependance r\'eguli\`ere en $\beta$ au
voisinage de l'origine. 
Pour une telle solution BKW \'el\'ementaire
$\Phi_{bkw}(z, \beta,\varepsilon)$, et pour $z\neq 0$ dans un voisinage
de $0$ et $\beta$ pr\`es de l'origine, 
 nous avons l'unique d\'ecomposition suivante~:
$$\Phi_{bkw}(z, \beta,\varepsilon) = a(z, \beta,\varepsilon) A_{bkw}(z,\varepsilon) +
b(z, \beta,\varepsilon)
\frac{\partial A_{bkw}}{\partial z}(z, \varepsilon),$$
o\`u $A_{bkw}(z,\varepsilon)$ est le symbole BKW d'Airy (\ref{wkb10}),
tandis que $a$ et $b$ sont des constantes locales de r\'esurgence.
\end{thm}

Nous traduisons \'egalement facilement le th\'eor\`eme \ref{theo42}.\\
Nous en d\'eduisons la cons\'equence suivante : en substituant \`a $\beta$ une
petite s\'erie Gevrey-1 $\displaystyle \beta (\varepsilon) = \sum_{k\geq 1}
\beta_k \varepsilon^{k}$ dans une constante locale de r\'esurgence, nous
obtenons encore une constante locale de r\'esurgence \cite{DP99}, nous retrouvons le
r\'esultat de Aoki {\em et al} \cite{Aoki93}.\\
De m\^eme, par extension, nous obtenons le th\'eor\`eme suivant~:
\begin{thm}\label{theo420}
Consid\'erons l'\'equation diff\'erentielle 
\begin{equation}\label{eqq5}
\displaystyle  \frac{d^2 \Phi}{d z^2} - \frac{z}{\varepsilon^2} \Phi = F(z, \varepsilon) \Phi,
\end{equation}
o\`u  $F(z, \varepsilon)$ est une constante locale de r\'esurgence . Alors
il existe une constante locale de r\'esurgence $s(z,\varepsilon)$ telle que,
sous l'action de la transformation
\begin{equation}\label{eq42ter}
 \left\{
 \begin{array}{l}
\displaystyle s(z, \varepsilon) = \sum_{k\geq 0} s_k(z)\varepsilon^{k}, \hspace{5mm} s_0(z) = z \\
\\
\displaystyle \Phi(z, \varepsilon)=\Big( \frac{\partial s}{\partial z} \Big)^{-\frac{1}{2}}
y \big( s(z, \varepsilon) , \varepsilon\big),
\end{array}
 \right.
\end{equation}
l'\'equation (\ref{eqq5}) est chang\'ee en l'\'equation (\ref{eq3}),
pour $z$ ({\em resp.} $s$) au voisinage de l'origine. De plus, sous
l'action de (\ref{eq42ter}), une solution BKW \'el\'ementaire de (\ref{eq3}) est
transform\'ee en une solution BKW \'el\'ementaire de (\ref{eqq5}).
\end{thm}

\section{Pistes de recherche}\label{sec6}

\subsection{Points tournants d'ordre sup\'erieur}

L'analogue de notre forme canonique pour les points tournants
d'ordre sup\'erieur est l'\'equation diff\'erentielle suivante~:
$$
\displaystyle  \frac{d^2 \Phi}{d z^2} - \frac{z^n}{\varepsilon} \Phi = F(z) \Phi,
\eqno{(M_n)}
$$
o\`u $F(z)$ est holomorphe au voisinage de l'origine et $n \in
\mathbb{N}\backslash \{0\}$.
Afin de copier ce que nous avons fait dans la section \ref{sec3}, il
nous faut d\'efinir une fonction g\'en\'eratrice convenable dans une
transformation canonique $(p,z) \leftrightarrow (\widehat p, \z)$ de
l'espace cotangent qui simplifie la g\'eom\'etrie de la
sous-vari\'et\'e Lagrangienne $p^2-z^n=0$. \\
Un point de d\'epart int\'eressant est l'article \cite{Hard10} o\`u
Hardy introduit un ensemble de fonctions sp\'eciales $\Phi_n$ solutions de~: 
$$
\displaystyle  \frac{d^2 \Phi}{d z^2} - \frac{z^n}{\varepsilon} \Phi = 0,
\eqno{(A_n)}
$$
sous la forme\footnote{
Hardy montre en particulier comment ces fonctions sont reli\'ees aux
fonctions de Bessel.}  
\begin{equation}\label{eq5ana}
\Phi_n (z, \varepsilon) =  \int  e^{-\frac{1}{\varepsilon} S_n (z, \z)} \, d\z.
\end{equation} 
Les fonctions $S_n$ peuvent \^etre d\'efinies de la mani\`ere suivante: pour $m=n+2 \geq 3$, nous
introduisons la fonction polynomiale $P_m(t) \in
\mathbb{Z}[t]$ d'ordre $m$ d\'efinie par~:
\begin{itemize}
\item  si  $m$ est pair, nous posons  $\cosh (mq) = P_m \left(\sinh
  (q)\right)$.
\item  si  $m$ est impair, nous posons  $\sinh (mq)=P_m \left( \sinh (q) \right)$.
\end{itemize}
Nous associons \`a $P_m$ la fonction polynomiale $Q_m$,
$$Q_m (z,\z)= z^{m/2}P_m(\z/\sqrt{z}),$$
et $S_n (z,\z) \in \mathbb{Q}[z,\z]$ est d\'efinie comme la fonction
polynomiale quasi-homog\`ene donn\'ee par~:
\begin{equation}\label{eq6ana}
S_n(z,\z)= \frac{2}{n+2}Q_{n+2} (-z,\z).
\end{equation}
Ces fonctions polynomiales $S_n$ satisfont les propri\'et\'es suivantes~:
\begin{equation}\label{eq7ana}
\left( \frac{\partial S_n}{\partial z} \right)^2 = 
T_n(z, \z)\frac{\partial S_n}{\partial \z } + z^n
\end{equation}
o\`u la fonction polynomiale $T_n (z, \z)$ satisfait~:
\begin{equation}\label{eq8ana}
\frac{\partial^2 S_n}{\partial z^2} = 
\frac{\partial  T_n}{\partial \z}.
\end{equation}
A titre d'exemple, nous avons~:
$$\left\{
\begin{array}{lll}
\vspace{2mm}
S_1(z, \z)=  \frac{8}{3}\z^3 -2\z z &  & T_1(z, \z)=  \frac{1}{2} \\
\vspace{2mm}
S_2(z, \z)= 4\z^4 -4\z^2 z + \frac{1}{2}z^2 & & T_2(z, \z)= \z \\
\vspace{2mm}
S_3(z, \z)=\frac{32}{5}\z^5 -8\z^3 z + 2\z z^2 & & T_3(z, \z)=2\z^2 -\frac{1}{2}z\\
 & \cdots &
\end{array}
\right.
$$
Le fait que la fonction $\Phi_n (z, \varepsilon)$ donn\'ee par (\ref{eq5ana})
soit en effet une solution de $(A_n)$ peut se montrer par int\'egration par parties,
en utilisant les propri\'et\'es fondamentales (\ref{eq7ana}) et (\ref{eq8ana}),  cf. \cite{C81}. \\
Notons que, \`a normalisation pr\`es,  $\Phi_1 (z,\varepsilon)$ correspond
\`a fonction d'Airy, tandis que $\Phi_2 (z,\varepsilon)$ correspond \`a la
fonction cylindro-parabolique de Weber.\\

En revenant \`a ce que nous avons fait dans la section \ref{sec3},
cela nous am\`ene \`a consid\'erer une solution de $(M_n)$ de la forme~:
\begin{equation}\label{eq9ana}
\Phi(z,\varepsilon) =  \int  e^{-\frac{\xi}{\varepsilon}} \stackrel{\vee}{\Phi}(z,\xi) \,
d\xi = \int  e^{-\frac{1}{\varepsilon} S_n(z, \z)} \Psi(z,\z)  \, d\z,
\end{equation}
o\`u  $\stackrel{\vee}{\Phi}$ doit satisfaire l'\'equation: 
\begin{equation}\label{eq10ana}
\frac{\partial^2 \stackrel{\vee}{\Phi}}{\partial z^2} - z^n \frac{\partial^2
  \stackrel{\vee}{\Phi}}{\partial \xi^2} = F(z) \stackrel{\vee}{\Phi}.
\end{equation}
En utilisant (\ref{eq7ana}) et (\ref{eq8ana}), ceci se traduit par le
fait que $\tpsi(z,\z) =
\stackrel{\vee}{\Phi} (z, \xi)$ doit satisfaire l'EDP lin\'eaire~:
\begin{equation}\label{eq11ana}
\frac{\partial^2 \tpsi}{\partial z^2} - \frac{2 \frac{\partial
    S_n}{\partial z} }{ \frac{\partial
    S_n}{\partial \z} }\frac{\partial^2
  \tpsi}{\partial z \partial \z} + \frac{T_n}{\frac{\partial
    S_n}{\partial \z} }\frac{\partial^2 \tpsi
 }{\partial \z^2} = F(z) \tpsi. 
\end{equation}
Comme dans la section \ref{sec3}, le probl\`eme se r\'eduit maintenant \`a analyser les
propri\'et\'es d'holomorphie de $\displaystyle  \Psi(z,\z) := \tpsi (z, \z)\frac{\partial
    S_n}{\partial \z}$ avec $\tpsi (z, \z)$ une solution de
  (\ref{eq11ana}) telle que  $\Psi(z,\z)$ se comporte bien au
  voisinage du lieu $\displaystyle \frac{\partial
    S_n}{\partial \z} =0$  d\'efinissant les points cols.\\
La difficult\'e r\'eside maintenant dans le fait que la transformation
$$  \left\{
\begin{array}{l}
\displaystyle (z,\z)\leftrightarrow(z,x=\frac{\partial S_n}{\partial
  \z})\\
\\
\displaystyle \psi(z,x):=\widetilde{\Psi}(z,\z)\frac{\partial
    S_n}{\partial \z},
\end{array} \right.$$
qui d\'efinit un changement de variables ramifi\'e, ne s'inverse plus
facilement (nous devons
r\'esoudre une \'equation alg\'ebrique de degr\'e $n+2$) et l'analyse
devient d'autant plus d\'elicate que $n$ devient grand. 

\subsection{Sommabilit\'e}

En ce qui concerne les propri\'et\'es de sommabilit\'e des solutions
BKW \'el\'ementaires du th\'eor\`eme \ref{propdecomp}, il semble, au regard des
exemples explicites obtenus dans la section \ref{sec412}, que nous perdions
assez vite le caract\`ere 1-sommable.\\
Notre conjecture est qu'en r\'ealit\'e nous obtenons des
solutions formelles multisommables, au moins dans le cas o\`u la
fonction $F$ est polynomiale (les travaux de Balser {\em et al} dans
\cite{} ne sont d'ailleurs certainement pas sans lien avec ce
ph\'enom\`ene).\\
Un travail est engag\'e dans cette pr\'esente voie afin de
d\'eterminer notamment les diff\'erents niveaux de sommabilit\'e (qui
doivent a priori d\'ependre du degr\'e de $F$).

\newpage

\setcounter{section}{0}
\renewcommand{\thesection}{\Alph{section}}

\section{Appendice : Fonctions confluentes et microfonctions}\label{sec7}

Nous rappelons ici quelques notions adapt\'ees de \cite{Jidoumou,
  DP99}.

\begin{defn}
\label{14}
Une microfonction  $\stackrel{\triangledown}{\Phi}(z, \xi)$ \`a $(z_0, \xi_0)$ est la classe 
modulo ${\mathcal O}_{
\mathbb{C} \times \mathbb{C} , z_0 \times \xi_0}$ d'une fonction 
${\stackrel{\vee}{\Phi}}(z, \xi)$  holomorphe dans un voisinage sectoriel
de  $(z_0, \xi_0)$. La fonction ${\stackrel{\vee}{\Phi}}$ est appel\'ee un majeur
et la microfonction correspondante sa singularit\'e \`a $(z_0,
\xi_0)$. \\

La microfonction $\stackrel{\triangledown}{\Phi}(z, \xi)$ est
r\'esurgente si son majeur  ${\stackrel{\vee}{\Phi}}$ est prolongeable
sans fin (par rapport \`a $\xi$ pour tout $z$ dans un voisinage de $z_0$).
\end{defn}

\begin{figure}[thp]
\begin{center}
\includegraphics[width=3.7in]{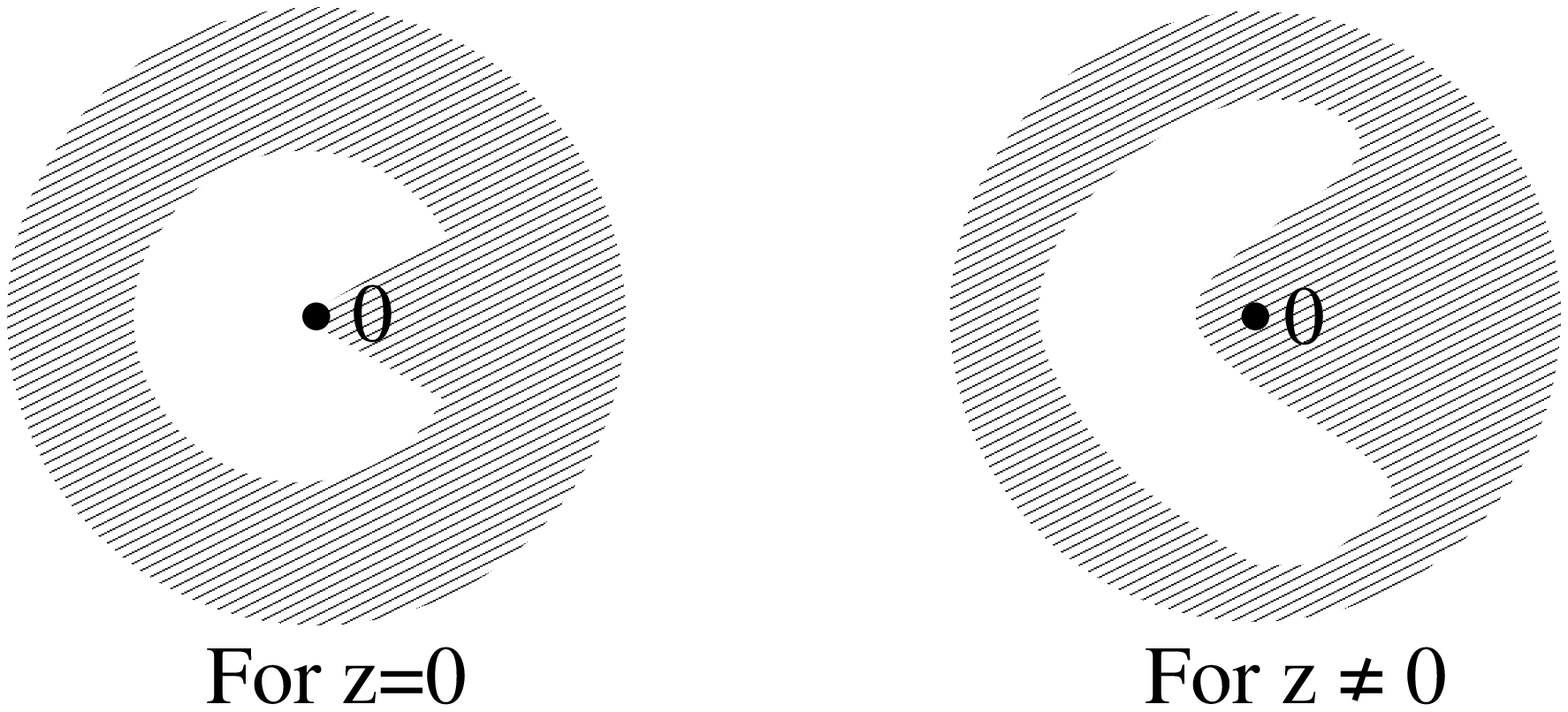} 
\caption{La trace dans le $\xi$-plan d'un voisinage sectoriel de $(0, 0)$.
\label{fig:Stokes7}}
\end{center}
\end{figure}

Dans cette d\'efinition, un  {\em voisinage sectoriel} de $(z_0,
\xi_0)$ d\'esigne un ouvert $W \subset \mathbb{C} \times \mathbb{C} $ intersectant 
$ \{z_0 \} \times \mathbb{C}$ comme le montre la figure \ref{fig:Stokes7}
pour $(z_0, \xi_0)=(0,0)$.

Nous rappelons que $\mathcal{C}$ d\'esigne la courbe alg\'ebrique $\mathcal{C} = \{(z, \xi) \in \mathbb{C}^2,\,\, 9\xi^2 = 4z^3 \}$. 
Nous allons utiliser la d\'efinition suivante~:

\begin{defn}
\label{15}
Une microfonction ( {\em  resp.}  microfonction r\'esurgente) $\stackrel{\triangledown}{\Phi}(z, \xi)$ 
\`a $(0, 0)$ est dite {\em confluente \`a support singulier dans
  $\mathcal{C}$} si l'un de ses majeurs ${\stackrel{\vee}{\Phi}}(z,
\xi)$ s'\'etend analytiquement sur le rev\^etement
universel de $\mathcal{U} \times \mathcal{V} \backslash
\mathcal{C}$ ({\em  resp.} de $\mathbb{C}^2 \backslash
\mathcal{C}$), o\`u
$\mathcal{U} \times \mathcal{V}$ est un voisinage de
$(0,0)$. Dans ce cas, ${\stackrel{\vee}{\Phi}}$ sera dit
{\em confluent} ( {\em  resp. confluent r\'esurgent}) \`a support
singulier dans $\mathcal{C}$. 
\end{defn}

Par exemple, la microfonction (r\'esurgente) \`a $(0, 0)$ associ\'ee \`a
$\displaystyle \stackrel{\vee}{A_{bkw}}(z,\xi)$ (cf. formule
(\ref{wkb11})) est confluente (r\'esurgente) \`a support singulier dans $\mathcal{C}$. 

\begin{figure}[thp]
\begin{center}
\includegraphics[width=1.8in]{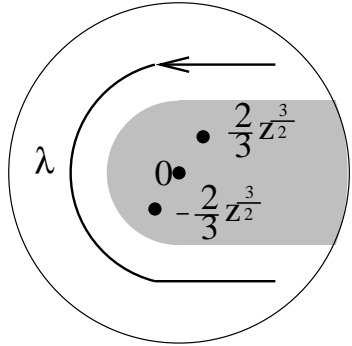} 
\caption{L'ouvert $\Omega$ (le compl\'ementaire de l'ensemble gris\'e)
  et le chemin $\lambda$. Les extr\'emit\'es de $\lambda$ doivent
  \^etre dans le demi-plan
  $\Re (\xi) >0$.
\label{fig:Stokes8}}
\end{center}
\end{figure}

Consid\'erons maintenant une microfonction confluente \`a $(0, 0)$ \`a
support singulier dans $\mathcal{C}$, et notons  ${\stackrel{\vee}{\Phi}}(z,
\xi)$ l'un de ses majeurs holomorphes dans un voisinage sectoriel de
$(0, 0)$. En gardant $z$ dans un voisinage suffisamment petit $\mathcal{U}$ de
$0$, nous pouvons d\'efinir un ouvert $\Omega$ comme sur la figure 
\ref{fig:Stokes8} tel que ${\stackrel{\vee}{\Phi}}$ est holomorphe
pour $(z,\xi) \in \mathcal{U} \times \Omega$.

Cela permet de d\'efinir la transform\'ee de Laplace suivante~:
\begin{equation}\label{eqa1}
\displaystyle \Phi(z,\varepsilon) =  \int_{\lambda}  e^{-\frac{\xi}{\varepsilon}}
\stackrel{\vee}{\Phi}(z,\xi) \, d\xi, 
\end{equation}
pour un contour d'int\'egration $\lambda$ comme sur la figure
\ref{fig:Stokes8}. La fonction $\displaystyle \Phi(z,\varepsilon)$ est
holomorphe pour $(z,\varepsilon) \in \mathcal{U} \times \mathbb{C}$, 
et admet une {\em croissance sous-exponentielle d'ordre 1 \`a l'infini en $\varepsilon^{-1}$} au
sens suivant~:\\
$\forall \eta >0$, il existe un voisinage
$\mathcal{U}_\eta$ de $0$, un voisinage sectoriel $\Sigma_\eta$ de l'origine d'ouverture $ ]-\delta_\eta,
\delta_\eta[$ ($\delta_\eta >0$) et $c_\eta >0$ tels que
\begin{equation}\label{eqa2}
\forall (z, \varepsilon) \in \mathcal{U}_\eta \times \Sigma_\eta, \, \, \,
\left| \Phi(z,\varepsilon) \right| \leq c e ^{\eta |\varepsilon|^{-1}}
\end{equation}

Si dans (\ref{eqa1}) nous rempla\c cons le majeur $\stackrel{\vee}{\Phi}$ par
par un autre repr\'esentant de la classe de la microfonction confluente
$\stackrel{\triangledown}{\Phi}$, ou changeons les extr\'emit\'es de
$\lambda$, la fonction $ \Phi(z,\varepsilon)$ est d\'ecal\'ee d'une
fonction holomorphe $\varphi(z, \varepsilon)$ en $z$ dans un voisinage $\mathcal{U}^\prime$ de $0$, et
en $\varepsilon \in \mathbb{C}^*$ qui est \`a {\em d\'ecroissance
  exponentielle d'ordre 1 \`a l'infini en $\varepsilon^{-1}$} au sens suivant~: \\
il existe $\eta >0$ et un voisinage sectoriel $\Sigma_\eta$ de l'origine d'ouverture $ ]-\delta_\eta,
\delta_\eta[$ ($\delta_\eta >0$), il existe $c>0$ tels que
\begin{equation}\label{eqa2bis}
\forall (z, \varepsilon) \in \mathcal{U}^\prime \times \Sigma_\eta, \, \, \,
\left| \varphi(z,\varepsilon) \right| \leq c e ^{-\eta |\varepsilon|^{-1}}
\end{equation}

Cela justifie la d\'efinition suivante \cite{Jidoumou}:

\begin{defn}\label{16}
La fonction confluente $\bold \Phi (z, \varepsilon)$ \`a support singulier dans 
$\mathcal{C}$ associ\'ee \`a la microfonction confluente $\stackrel{\triangledown}{\Phi}(z, \xi)$ 
\`a $(0, 0)$ est la classe de $\Phi (z, \varepsilon)$ d\'efinie par (\ref{eqa1})
modulo les fonctions holomorphes en $(0,0)$ qui sont \`a
d\'ecroissance exponentielle d'ordre 1 \`a l'infini en $\varepsilon^{-1}$.
\end{defn}

\begin{rem}\label{16bis}
Si ${\stackrel{\vee}{\Phi}}$ est confluent r\'esurgent \`a support singulier
dans $\mathcal{C}$, alors la fonction confluente r\'esurgente $\bold \Phi
(z, \varepsilon)$ \`a support singulier dans $\mathcal{C}$ associ\'ee \`a la
microfonction confluente r\'esurgente $\stackrel{\triangledown}{\Phi}(z, \xi)$ 
\`a $(0, 0)$ est d\'efinie comme une pr\'esomme de Borel, voir \cite{DP99, CNP2}.
\end{rem}

\subsection*{A.1 : D\'ecomposition locale (pour la direction $\alpha =0$)} 
\addcontentsline{toc}{subsubsection}{A.1 : D\'ecomposition locale (pour la
  direction $\alpha=0$)}

Pour la notion de d\'ecomposition d'une fonction r\'esurgente, nous
renvoyons le lecteur \`a \cite{DP99, CNP2}.
Nous introduisons ici son analogue local.
Consid\'erons de nouveau un majeur  ${\stackrel{\vee}{\Phi}}(z,
\xi)$ dans la classe d'une microfonction confluente \`a $(0, 0)$ \`a
support singulier dans $\mathcal{C}$. Nous gardons $z$ dans un voisinage
suffisamment petit $\mathcal{U}$ de
$0$ et $z \neq 0$. Pour un tel $z$, la fonction 
$\xi \mapsto {\stackrel{\vee}{\Phi}}(z,
\xi)$ peut se prolonger dans un domaine coup\'e (localement pr\`es de $\xi
=0$) comme sur la figure \ref{fig:Stokes5}.b,  (cas I), ou sur la
figure \ref{fig:Stokes62} (cas II). \\
En consid\'erant les prolongements analytiques
de $\xi \mapsto {\stackrel{\vee}{\Phi}}(z,
\xi)$ au voisinage de chaque singularit\'e $\omega (z)$ 
(d\'efinie par $\mathcal{C}$), et les microfonctions associ\'ees (i.e., la classe 
modulo $O_{\omega(z)}$, $\omega (z)$ \'etant le {\em support} de la microfonction), nous
obtenons la {\em d\'ecomposition locale } d'une microfonction confluente
$\stackrel{\triangledown}{\Phi}$. Par exemple, dans le cas I, la
d\'ecomposition locale est donn\'ee par seulement une microfonction 
$\stackrel{\triangledown}{\Phi}^{\frac{2}{3}z^{3/2}}$
\`a $\frac{2}{3}z^{3/2}$, alors que dans le cas II, la d\'ecomposition locale 
fait intervenir deux microfonctions, l'une  $\stackrel{\triangledown}{\Phi}^{\frac{2}{3}z^{3/2}}$
\`a $\frac{2}{3}z^{3/2}$, et l'autre $\stackrel{\triangledown}{\Phi}^{-\frac{2}{3}z^{3/2}}$
\`a $-\frac{2}{3}z^{3/2}$.

Nous utilisons maintenant les germes de secteurs de Stokes, comme sur la 
 figure \ref{fig:Stokes5}.a. Tant que nous restons dans un tel
 germe de secteurs de Stokes $S$, les microfonctions
 $\stackrel{\triangledown}{\Phi}^{\omega(z)}$ intervenant dans la d\'ecomposition locale
  de la microfonction confluente
$\stackrel{\triangledown}{\Phi}$ sont holomorphes en $z$. Cela nous
permet de d\'efinir le morphisme
\begin{equation}\label{eqa3}
\begin{array}{ccc} 
\stackrel{\triangledown}{\Phi}(z, \xi) &
\stackrel{\displaystyle \sigma_S}{\displaystyle \longrightarrow} &
\left(\stackrel{\triangledown}{\Phi}^{\omega(z)}(z, \xi) \right)_{\omega(z)}
\end{array}
\end{equation}
A chaque  microfonction $\stackrel{\triangledown}{\Phi}^{\omega(z)}(z,
\xi)$ nous pouvons associer, par l'interm\'ediaire d'une transform\'ee
de Laplace formelle, un unique symbole BKW \'el\'ementaire
$\Phi_{bkw}^{\omega(z)}(z, \varepsilon)$. 
Par suite, nous avons le diagramme commutatif suivant~:
\begin{equation}\label{eqa4}
\begin{array}{ccc}
\stackrel{\triangledown}{\Phi}(z, \xi) &
\stackrel{\displaystyle \sigma_S}{\displaystyle \longrightarrow} &
\left(\stackrel{\triangledown}{\Phi}^{\omega(z)}(z, \xi) \right)_{\omega(z)}\\
\downarrow & & \downarrow\\
\bold \Phi (z, \varepsilon) & \stackrel{\displaystyle \sigma_S}{\displaystyle
  \longrightarrow} & \sum_{\omega(z)} \Phi_{bkw}^{\omega(z)}(z, \varepsilon)
\end{array}
\end{equation}
o\`u $\displaystyle \sum_{\omega(z)} \Phi_{bkw}^{\omega(z)}(z,
\varepsilon)$, qui d\'esigne la s\'erie transasymptotique locale
associ\'ee \`a la fonction confluente $\bold \Phi (z, \varepsilon)$, est
appel\'ee la {\em d\'ecomposition locale
 dans $S$ de la fonction confluente $\bold \Phi (z,
\varepsilon)$}.

Lorque nous traversons une ligne de Stokes, un ph\'enom\`ene de Stokes
appara\^it. Ce dernier se traduit par une discontinuit\'e de la
d\'ecomposition locale, et peut s'analyser en termes de d\'erivations \'etrang\`eres.
Cette \'etude a d\'ej\`a \'et\'e men\'ee par exemple pour le cas Airy
dans la section \ref{sec2}.


\end{document}